%% file: skein-dimensions.tex
\documentclass[a4paper]{amsart}
\usepackage[utf8]{inputenc}
\usepackage[left=2.5cm, right = 2.5cm, top = 3.5cm, bottom=3.5cm]{geometry}

\usepackage[dvipsnames]{xcolor}
\definecolor{darkblue}{rgb}{0, 0, 0.6}
\usepackage[pdftitle={\@title}, pdfauthor={\@author}, pdfpagelayout=OneColumn, ocgcolorlinks, colorlinks=true, linkcolor=darkblue, urlcolor=NavyBlue, citecolor=darkblue, filecolor=darkblue]{hyperref}
\usepackage[style=alphabetic, citestyle=alphabetic, maxnames=4, maxalphanames=4]{biblatex}
\DeclareBibliographyCategory{needsurl}
\renewbibmacro*{url+urldate}{%
  \ifcategory{needsurl}
    {\printfield{url}%
    }
    {}}
\newcommand{\entryneedsurl}[1]{\addtocategory{needsurl}{#1}}
\entryneedsurl{walkerTQFTs} 
\entryneedsurl{Kin25SkeinModulesTwisted}
\entryneedsurl{hatcherNotesBasic3Manifold}

\AtEveryBibitem{\clearfield{urlyear}}

\DeclareFieldFormat{url}{\newline\mkbibacro{URL}\addcolon\nobreakspace\url{#1}}

\usepackage{import}
\usepackage[toc]{appendix}

\usepackage{graphicx}
\usepackage{subcaption}
\usepackage{float}
\usepackage{listings}
\usepackage{stmaryrd}
\usepackage{pkinne}

\graphicspath{{images}}
\usepackage[notransparent, inkscapepath=svgsubdir]{svg}

\title{Skein module dimensions of mapping tori of the 2-torus}

\author{Patrick Kinnear}
\address{School of Mathematics, University of Edinburgh, Edinburgh, UK}
\email{P.Kinnear@ed.ac.uk; patrick.kinnear@uni-hamburg.de}


\newcommand\SkAlg{\mathrm{SkAlg}}
\newcommand\SkCat{\mathrm{SkCat}}
\newcommand\Sk{\mathrm{Sk}}

\def\fd{\mathrm{f.d.}}
\def\Stab{C_W}
\def\HH{\mathrm{HH}}
\def\H{\mathrm{H}}
\def\tr{\mathrm{tr}}
\def\rk{\mathrm{rk}}

\newcommand{\conj}[1]{\mathrm{cl}({#1})}

\def\Kar{\mathrm{Kar}}
\def\Free{\mathrm{Free}}
\def\TL{\mathcal{TL}}

\numberwithin{equation}{section}

\begin{document}

\begin{abstract}
  We determine the dimension of the Kauffman bracket skein module at generic $q$ for mapping tori of the 2-torus, generalizing the well-known computation of Carrega and Gilmer. In the process, we give a decomposition of the twisted Hochschild homology of the $G$-skein algebra for $G = \SL_N$ or $\GL_N$, which is a direct summand of the whole skein module, and from which the dimensions follow easily in the cases $G = \SL_2$ and $G = \GL_1$.
\end{abstract}

\maketitle

\section{Introduction}
  \import{}{1-intro/outline}
  \subsection*{Acknowledgements}
  \import{}{1-intro/ack}

\section{Background and preliminaries}
  \label{s-prelim-skeins}
  \import{}{2-prelim/outline}
  \subsection{Mapping tori of the 2-torus}
  \label{s-mapping-tori}
  \import{}{2-prelim/mapping-tori}
  \subsection{Twisted Hochschild homology of skein categories}
  \label{s-twisted-HH_0}
  \import{}{2-prelim/twisted-HH_0}
  \subsection{Representation-theoretic description of the skein algebra}
  \label{s-prelim-SkAlg}
  \import{}{2-prelim/SkAlg}

\section{Twisted Hochschild homology of the skein algebra}
  \label{s-HH_0-SkAlg}
  \import{}{3-HH_0-SkAlg/outline}
  \subsection{Direct sum decomposition}
  \label{s-HH_0_decomp}
  \import{}{3-HH_0-SkAlg/decomp}
  \subsection{Twisted Hochschild homology in terms of the difference cokernel}
  \label{s-diff-cokernel}
  \import{}{3-HH_0-SkAlg/diff-cokernel}

\section{Calculations for the Kauffman bracket skein module}
  \label{s-SL_2-dims}
  \import{}{4-SL_2-dims/outline}
  \subsection{Single skein part}
  \label{s-single-skein}
  \import{}{4-SL_2-dims/single-skein}
  \subsection{Total skein module dimensions}
  \label{s-SL_2-dims-total}
  \import{}{4-SL_2-dims/total}

\clearpage
\begin{appendices}
  \section{Sage implementation}
  \label{a-sage}
  \import{}{A-Appendix/code}
  \section{Table of dimensions}
  \label{a-table}
  \import{}{A-Appendix/table}
  \section{A simpler formula}
  \label{a-simpler}
  \import{}{A-Appendix/derivation}
\end{appendices}

\sloppy
\printbibliography

\end{document}

%% file: 1-intro/outline.tex
\begin{figure}[h]
  \begin{center}
    \includesvg[width=7cm]{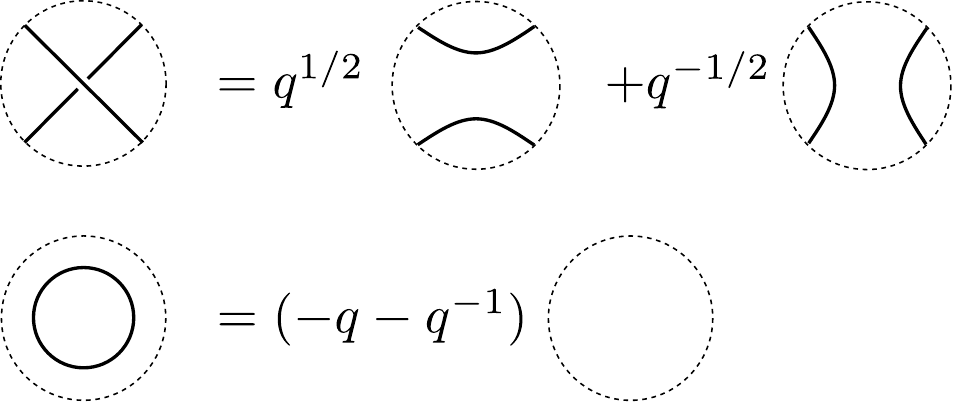}
    \caption{The Kauffman bracket skein relations.}
    \label{f-KBSRs}
  \end{center}
\end{figure}

The Kauffman bracket (or $\SL_2$) skein module $\Sk_{\SL_2}(M)$ of a 3-manifold $M$ is the $\C(q^{1/2})$-vector space spanned by isotopy
classes of framed links embedded in $M$, modulo the skein relations depicted in Fig. \ref{f-KBSRs}.
It was conjectured by Witten that skein modules of closed 3-manifolds should
be finite-dimensional over $\C(q^{1/2})$. This question was only recently settled by Gunningham--Jordan--Safronov \cite{gunninghamFinitenessConjectureSkein2019},
however the proof does not give access to explicit dimensions, which have still only been computed in a
few cases. The dimensions have been computed for $S^3$ and lens spaces \cite{hosteSKEINMODULELENS1993, hosteKauffmanBracketSkein1995}, integer Dehn surgeries along a trefoil \cite{bullockKauffmanBracketSkein1997}, the quaternionic manifold \cite{gilmerKauffmanBracketSkein2007}, some prism manifolds \cite{mroczkowskiKAUFFMANBRACKETSKEIN2011}, and some families of hyperbolic manifolds obtained by Dehn filling of knot complements \cite{DKS23KauffmanBracketSkein}. The dimension was computed for the 3-torus by Carrega \cite{carregaGeneratorsSkeinSpace2017} and Gilmer \cite{gilmerKauffmanBracketSkein2018}, and this was generalized to the product of a closed surface with a circle in work of Gilmer--Masbaum \cite{gilmerSkeinModuleProduct2019a} and Detcherry--Wolff \cite{detcherryBasisKauffmanSkein2021}. In the current work, we make a contribution to this growing list of 3-manifolds by computing the Kauffman bracket skein module dimension for the infinite family of 3-manifolds which are mapping tori $M_\gamma = T^2 \times_\gamma S^1$ of mapping classes $\gamma \in \mathrm{Mod}(T^2)$ of the 2-torus $T^2$.

Recalling the abelian case, the $\GL_1$-skein module of a manifold $M$ is already well-understood due to work of Przytycki \cite{Prz98AnalogueFirstHomology}: it is isomorphic to the $\C(q)$-vector space supported on the torsion part of $\H_1(M)$, which for mapping tori $M_\gamma$ is easily computed (see \S \ref{s-mapping-tori}). Our main result can be seen as a generalization of this calculation to the non-abelian setting. We also prove some supporting results which allow us to decompose the $\GL_N$- and $\SL_N$-skein modules of mapping tori: in particular we compute the $\GL_1$-skein module dimensions and show that they recover the calculations of Przytycki (see \S \ref{s-diff-cokernel}). We note that the $\SL_N$ and $\GL_N$-skein module dimensions for the 3-torus were recently computed in \cite{gunninghamSkeinsTori}, and viewing $T^3$ as $M_{\Id}$, our work can be regarded as generalizing the techniques of \cite{gunninghamSkeinsTori} for nontrivial mapping tori.

\subsection{Main results}

\begin{notn}
  Let $\gamma \in \Mod(T^2) \cong \SL_2(\Z)$, and denote by $\bar{\gamma}$ the matrix for $\gamma$ taken modulo 2 in each entry. Denote by $a^{\pm}_i$ the invariant factors of the map $\Id \mp \gamma : \H_1(T^2) \to \H_1(T^2)$ and $r_{\pm}$ the rank of this map. Let $p_{\pm} = \# \{a_i^{\pm} \mathrm{ even} : 1 \leq i \leq r_{\pm} \}$.
\end{notn}

\begin{thm}
  \label{th-intro-dimension-formula}
  The dimension of the Kauffman bracket skein module of $M_\gamma$ is given by is given by
  \[
    \dim \Sk_{\SL_2}(M_\gamma) = s_\gamma + \frac{\prod_{i = 1}^{r_+}a_i^{+} + 2^{p_{+}}}{2} + \frac{\prod_{i = 1}^{r_-}a_i^{-} + 2^{p_{-}}}{2} 
  \]
  where
  \[
    s_\gamma = \begin{cases}
      4 & \bar{\gamma} = \bar{\Id}\\
      2 & \tr(\bar{\gamma}) = 0 \mod 2 \text{ and } \bar{\gamma} \neq \bar{\Id}\\
      1 & \tr(\bar{\gamma}) = 1 \mod 2.
    \end{cases}
  \]
\end{thm}

\begin{rmk*}
  After the publication of this paper, the author realized that the above formula admits a significant simplification when $|\tr(\gamma)| > 2$. In this case, the formula can be written as
  \[
    \dim \Sk_{\SL_2}(M_{\gamma}) = |\tr(\gamma)| + 2^{c(\gamma) + 1}
  \]
  where
  \[
    c(\gamma) = \#\{m \in \{\mathrm{gcd}(a-1, b, c, d-1), \tr(\gamma)\} : m \text{ even} \}
  \]
  and $a, b, c, d$ are the entries of $\gamma$. The above formula can be derived in an elementary way from Theorem \ref{th-intro-dimension-formula}. This simple formula renders the large table of dimensions and computer implementation included in the paper unnecessary. However, we include these and the rest of the paper unedited for continuity with the published version, and add a further appendix (Appendix \ref{a-simpler}) with a short derivation of the above simple formula (which first appeared in \cite{Kin25SkeinModulesTwisted}).
\end{rmk*}

\begin{eg}
  \label{eg-SL_2-formulas}
  The formula of Thm. \ref{th-intro-dimension-formula} is straightforward to implement on a computer (see Appendix \ref{a-sage} for a Sage implementation).
  We recover the well-known theorem of \cite{carregaGeneratorsSkeinSpace2017, gilmerKauffmanBracketSkein2018} that $\dim \Sk_{\SL_2}(T^3) = 9$. As described in \S \ref{s-mapping-tori}, mapping tori are classified up to oriented diffeomorphism by conjugacy classes in $\SL_2(\Z)$, and as noted in Rmk. \ref{r-ccls-up-to-sign} the skein modules of $M_{\pm \gamma}$ have the same $\SL_2$-skein module dimensions.  For conjugacy classes $[\gamma]$ with $0 \leq \tr(\gamma) \leq 2$, the dimensions are tabulated in Table \ref{t-lowtrace-dims}. For mapping classes of trace greater than 2, we tabulate some selected skein dimensions in Table \ref{t-length-2-dims}.
\end{eg}

\begin{table}
  \begin{center}
    \begin{tabular}{c | c || c c c | c}
      $\tr(\gamma)$ & $\gamma$                                  & $s_\gamma$ & $d_{+}$ & $d_{-}$ & $\dim \Sk_{\SL_2}(M_\gamma)$ \\
      \hline
      0             & $S = \begin{pmatrix}0 & -1\\ 1 & 0\end{pmatrix}$           & 2          & 2             & 2             & 6                            \\
      1             & $E_+ = \begin{pmatrix}1 & -1\\ 1 & 0\end{pmatrix}$         & 1          & 1             & 2             & 4                            \\
      1             & $E_- = \begin{pmatrix}0 & 1\\ -1 & 1\end{pmatrix}$         & 1          & 1             & 2             & 4                            \\
      2             & $T^{2k} = \begin{pmatrix}1 & 2k\\ 0 & 1\end{pmatrix}$      & 4          & $k + 1$       & 4             & $9 + k$                      \\
      2             & $T^{2k + 1} = \begin{pmatrix}1 & 2k + 1\\ 0 & 1\end{pmatrix}$ & 2          & $k + 1$       & 3             & $6 + k$                      \\
    \end{tabular}
  \end{center}
  \caption{Kauffman bracket skein module dimensions for low-trace mapping classes: as explained in \S \ref{s-mapping-tori}, these represent all mapping tori $M_\gamma$ with $0 \leq \tr(\gamma) \leq 2$, up to oriented diffeomorphism. The two summands depending on the invariant factors of $\Id \mp \gamma$ are denoted $d_{\pm}$. As noted in Rmk. \ref{r-ccls-up-to-sign}, the $\SL_2$ dimensions for $-\gamma$ will be the same, with the roles of $d_{\pm}$ reversed. Note the dimension for $\gamma = \Id = T^{0}$ agrees with the known result of \cite{carregaGeneratorsSkeinSpace2017, gilmerKauffmanBracketSkein2018} for the 3-torus.}
  \label{t-lowtrace-dims}
\end{table}

Our calculation of $\SL_2$-skein dimensions proceeds by giving a direct sum decomposition of $\Sk_G(M)$ in the more general case $G = \GL_N$ or  $\SL_N$. Each of the direct summands in this decomposition can be interpreted as a $\gamma$-twisted version of the Hochschild homology of some algebra. One of these algebras will always be the $G$-skein algebra of $T^2$ (let us note that the Hochschild homology of the skein algebra has been considered in \cite{oblomkovDoubleAffineHecke2004,mclendonDetectingTorsionSkein2006, mclendonTracesSkeinAlgebra2007}). Our first main theorem is a decomposition of the $\gamma$-twisted Hochschild homology of the $G$-skein algebra.

\begin{thm}
  \label{th-intro-HH_0-decomp}
  Let $G = \GL_N$ or $\SL_N$, and $\gamma \in \Mod(T^2)$. Then
  \[
    \HH_0^\gamma(\SkAlg_G(T^2)) \cong \bigoplus_{[w] \in \conj{W}} \HH_0(A, A_{\gamma w})_{\Stab(w)}
  \]
  where $W$ is the Weyl group of $G$ and $A$ is a quantum torus algebra (see Notation \ref{n-quantum-torus}), and $A_{\gamma w}$ a twisted bimodule (see Def. \ref{d-twisted-HH_0}). Here, the direct sum ranges over the set $\conj{W}$ of conjugacy classes in the Weyl group, and the subscript $\Stab(w)$ denotes the space of coinvariants for the centralizer subgroup of an element $w \in W$.
\end{thm}

This allows us to give a general decomposition for skein modules of mapping tori for $\GL_N$ and $\SL_N$.  In the $\GL_N$ case only the twisted Hochschild homology  of the skein algebra is required; in the $\SL_N$ case some further endomorphism algebras in the skein category appear.

\begin{cor}
  \label{c-intro-skmod-decomp}
  For $G = \GL_N$, we have a decomposition of the skein module of $M_\gamma$ as
  \[
    \Sk_{\GL_N}(M_\gamma) \cong \HH_0^\gamma(\SkAlg_{\GL_N}(M)) \cong \bigoplus_{[w] \in \conj{W}} \HH_0(A, A_{\gamma w})_{\Stab(w)}.
  \]
  For $G = \SL_N$, we have the decomposition
  \[
    \Sk_{\SL_N}(M_\gamma) \cong \bigoplus_{n = 1}^{N-1} \HH_0^{\gamma}(\End_{\SkCat_{\SL_N}(T^2)}(V^{\ot n})) \oplus \bigoplus_{[w] \in \conj{W}} \HH_0(A, A_{\gamma w})_{\Stab(w)}
  \]
  where $V^{\ot n}$ denotes the object of the skein category of $T^2$ given by $n$ marked points on $T^2$  each coloured by the defining representation $V$ of $U_q(\sl_N$) (see \S \ref{s-twisted-HH_0}).
\end{cor}

In the case $G = \SL_2$, it is straightforward to compute the spaces of coinvariants appearing in our decomposition, which deals with the twisted Hochschild homology of the skein algebra, and gives the contributions $n_i^{\pm}$ in Thm. \ref{th-intro-dimension-formula}. In this case there is only one other direct summand of the skein module which gives the contribution $s_\gamma$ in Thm. \ref{th-intro-dimension-formula}.

\begin{rmk}
  \label{r-trace-property}
  Skein theory has been situated in the framework of a 3-2-TQFT \cite{walkerTQFTs,johnson-freydHeisenbergpictureQuantumField2019}. Our theorem gives an easy proof that skein theory cannot be extended to an oriented 4-3-2-TQFT such that the traces of the associated representation of $\Mod(T^3)$ are bounded above. Suppose there was such a TQFT assigning to a 3-manifold its skein module. Then the partition function of a mapping torus of $T^3$ (a 4-manifold) should be the trace of the mapping class of $T^3$ acting on $\Sk_{\SL_2}(T^3)$. Considering mapping classes $\gamma$ of $T^2$ under the embedding $\SL_2(\Z) \hookrightarrow \SL_3(\Z)$, these give mapping classes of $T^3$ and the corresponding mapping torus has the form $(T^2\times_{\gamma} S^1) \times S^1 \cong M_{\gamma} \times S^1$. The partition function of this 4-manifold should be the dimension of $\Sk_{\SL_2}(M_{\gamma})$. We can observe from the computations presented in Tables \ref{t-lowtrace-dims} and \ref{t-length-2-dims} that these dimensions are not bounded above, which contradicts the assumption that the 4-3-2-TQFT gives a representation of $\Mod(T^3)$ with bounded trace. So no such TQFT exists.
  
  There is a natural action of $\Mod(T^3) \cong \SL_3(\Z)$ on the skein module of $T^3$, where a mapping class acts on a skein in $T^3$ by acting on $T^3$. We recall the observation \cite[Rmk. 2.17]{carregaGeneratorsSkeinSpace2017} that this action is given by permutations of the 9 basis elements. Then the trace of this action will be an integer less than or equal to 9, so the natural action of $\Mod(T^3)$ has bounded trace and cannot be obtained from a 4-3-2-TQFT which extends skein theory. This implication was first pointed out to us by R. Detcherry.
\end{rmk}

\begin{rmk}
  The invariant factors $a_i^{\pm}$ appearing in Thm. \ref{th-intro-dimension-formula} give the direct sum decomposition of the finitely generated abelian group $\coker(\Id \mp \gamma)$, which in \cite[\S 3.2.1]{CGPS203d3dCorrespondenceMapping} is identified with the set of abelian ($\Id - \gamma$) and almost abelian ($\Id + \gamma$) flat $\SL_2(\C)$-connections on $M_{\gamma}$. It would be interesting to understand more fully how Thm. \ref{th-intro-dimension-formula} relates to calculations inspired by the 3d-3d correspondence of \cite{DGH11VortexCountingLagrangian}, such as the so-called $\hat{Z}$-invariants which were studied for mapping tori in \cite{CGPS203d3dCorrespondenceMapping}. As in Rmk. \ref{r-trace-property} this may reveal something about the 4d TQFT categorifying $\hat{Z}$-invariants.
\end{rmk}

\begin{qu}
  We emphasize that Cor. \ref{c-intro-skmod-decomp} gives a decomposition of $\Sk_G(M)$ for any $G = \GL_N$ or $\SL_N$. In the case of $\SL_2$ the centralizers of the Weyl group which appear are straightforward to handle. Gaining an understanding of the combinatorics of the centralizers in the general case would yield the $\GL_N$ dimensions and the summand of the $\SL_N$ corresponding to the skein algebra. In the $\SL_N$ case, the total dimension will depend on $N-1$ further summands, which will not in general be as easy to handle as the $\SL_2$ case. In particular, as identified in \cite{gunninghamSkeinsTori}, for composite $N$ there will be summands which are the twisted Hochschild homology of an infinite dimensional algebra. We defer these problems to future researchers.
\end{qu}

\subsection{Layout of the paper}

In \S \ref{s-prelim-skeins}, after recalling some background on mapping tori, we explain how to decompose $\Sk_G(M_\gamma)$ as a direct sum of twisted Hochschild homologies of endomorphism algebras in the skein category, building on work of \cite{gunninghamSkeinsTori}. We take care to set up the definitions of skein categories, modules and algebras, in the case of a general group $G$. We moreover introduce the algebra $A$ that appears in the above statements, and is such that $A^W$ is the $G$-skein algebra of $T^2$.

In \S \ref{s-HH_0_decomp} we show how to decompose the direct summand corresponding to the skein algebra giving Thm. \ref{th-intro-HH_0-decomp}. We explain in \S \ref{s-diff-cokernel} how to change basis in $A$ to understand the spaces $\HH_0(A, A_{\gamma w})$. In the case $G = \GL_1$, the Weyl group is trivial and it is immediate to obtain the $\GL_1$-skein module dimensions.

In the case $G = \SL_2$, it is straightforward to compute the spaces of coinvariants appearing in our decomposition, which handles the summand corresponding to the skein algebra. The remaining direct summand is dealt with in \S \ref{s-single-skein}. In \S \ref{s-SL_2-dims-total} we collect our computations to give the proof of Thm. \ref{th-intro-dimension-formula}.

%% file: 1-intro/ack.tex
The author would like to thank his advisors: David Jordan for his support and guidance throughout the project, and Pavel Safronov for many helpful discussions. The author also thanks Renaud Detcherry for his interest in the project and valuable conversations and feedback, Iordanis Romaidis for his suggestions, and Alisa Sheinkman for her early involvement. The author is grateful to the anonymous reviewers, whose comments have improved the paper. The author was supported by the Carnegie Trust for the Universities of Scotland for the duration of this research.

%% file: 2-prelim/outline.tex
In this section we recall some useful background. In \S \ref{s-mapping-tori} we give some details on mapping tori of $T^2$ that are not essential to our main computations but may serve as helpful orientation. In \S \ref{s-twisted-HH_0} we introduce the notion of twisted Hochschild homology of the skein category: we explain how this is used to obtain the skein module of a mapping torus, and how it decomposes as a direct sum. One direct summand is the twisted Hochschild homology of the skein algebra of $T^2$, and in \S \ref{s-prelim-SkAlg} we recall a useful description of this skein algebra. In \S \ref{s-HH_0-SkAlg} we will give a decomposition of this twisted Hochschild homology which holds for the $G$-skein algebra, when $G = \SL_N$ or $\GL_N$. We are therefore careful in this section to introduce notions such as skein categories and skein algebras for general groups $G$.

When $G$ does not have even Cartan determinant, there is some technical care required in treating the ground field, which we fix here.

\begin{notn}
  \label{n-fields}
  We work over an algebraically closed field $k$ of characteristic zero. When we consider the $G$-skein module $\Sk_G(M)$, this is a vector space naturally defined over $k = \C(q^{1/d})$ where $d$ is the determinant of the symmetrized Cartan matrix of $G$. In speaking of the skein module dimension for generic $q$, we mean the dimension of $\Sk_G(M)$ over $k$. In places we need to assume that $k$ contains the element $q^{1/2}$. Where $2$ does not divide $d$, then we work with $\cK = \C(q^{1/d}, q^{1/2})$ and the vector space $\cK \ot_k \Sk_G(M)$. Determining the $\cK$-dimension of $\cK \ot_k \Sk_G(M)$ is equivalent to determining the $k$-dimension of $\Sk_G(M)$. We will typically suppress mention of the field when discussing dimensions, on the grounds that the dimension will always be computed working over the field which makes technical sense and this dimension will agree with the dimension of the vector space over the field for which it is naturally defined. In the main case of interest where $G = \SL_2(\C)$, we have $d =2$ and $k = \cK = \C(q^{1/2})$. More generally, $d = 1$ for $\GL_N$ and $d = N$ for $\SL_N$.
\end{notn}

%% file: 2-prelim/mapping-tori.tex
We are interested in mapping tori of $T^2$, which we denote by $M_\gamma$ where $\gamma$ is the mapping class in question. These manifolds are solvmanifolds when the monodromy $\gamma$ is Anosov, and are known Seifert manifolds otherwise. Here we recall the classification of mapping tori, giving details on their Seifert description where it exists. Moreover we recall the computation of the first homology groups of these mapping tori.

\begin{defn}
  Given a mapping class $\gamma \in \Mod(T^2)$, the mapping torus of $\gamma$ is defined as
  \[
    M_\gamma = (T^2 \times I) / ((a, 0) \sim (\gamma(a), 1)).
  \]
\end{defn}

Recall that $\Mod(T^2) \cong \SL_2(\Z) \subseteq \GL_2(\Z)$ and $\SL_2(\Z)$ has a presentation
\[
  \SL_2(\Z) = \langle S, T | S^4 = 1, (ST)^3 = S^2 \rangle
\]
given by taking $S = \begin{pmatrix}0 & -1\\ 1 & 0\end{pmatrix}, T = \begin{pmatrix}1 & 1\\ 0 & 1\end{pmatrix}$. Two mapping tori $M_\gamma, M_\phi$ are diffeomorphic as oriented manifolds if and only if the monodromies $\gamma, \phi$ are conjugate in $\SL_2(\Z)$; they are diffeomorphic (possibly reversing orientation) if and only if $\gamma$ is conjugate to $\phi^{\pm}$ in $\GL_2(\Z)$ \cite[Thm. 2.6]{hatcherNotesBasic3Manifold}. Therefore oriented diffeomorphism classes of mapping tori correspond to conjugacy classes in $\SL_2(\Z)$.

The conjugacy classes in $\SL_2(\Z)$ are classified as follows (see, e.g. \cite[Chapter 7]{karpenkovContinuedFractionsOperatornameSL2013}), where we write $\gamma \sim \phi$ to denote the conjugacy relation.
\begin{itemize}
  \item $\tr(\gamma) = 0$: then $\gamma \sim \pm S$. The mapping torus of $S$ is diffeomorphic to the Seifert manifold $\{-2; (o_1, 0); (2, 1), (4, 3), (4, 3)\}$ \cite[\S 8.2]{orlikSeifertManifolds1972}.
  \item $\tr(\gamma) = 1$: then $\gamma \sim E_+ = \begin{pmatrix}1 & -1\\ 1 & 0\end{pmatrix}$ or $\gamma \sim E_- = \begin{pmatrix}0 & 1\\ -1 & 1\end{pmatrix}$, which have order 6.
  The corresponding mapping tori are Seifert manifolds: the mapping torus for $E_+$ is diffeomorphic to $\{-2; (o_1, 0); (2, 1), (3, 2), (6, 5) \}$, and the mapping torus for $E_-$ is $\{-1; (o_1, 0); (2, 1), (3, 1), (6, 1) \}$ \cite[\S 8.2]{orlikSeifertManifolds1972}.
  \item $\tr(\gamma) = -1$: then $\gamma \sim -E_+$ or $\gamma \sim -E_-$. These have order 3. The corresponding manifolds are Seifert manifolds: the mapping torus for $-E_+$ is $\{-1; (o_1, 0); (3, 1), (3, 1), (3, 1) \}$, and the mapping torus for $-E_-$ is $\{-2; (o_1, 0); (3, 2), (3, 2), (3, 2) \}$ \cite[\S 8.2]{orlikSeifertManifolds1972}.
  \item $\tr(\gamma) = 2$: then there is a $\Z$-indexed family of conjugacy classes given by $\{ T^n : n \in \Z \}$. We call these equivalence classes of mapping classes the shears. The corresponding mapping tori are diffeomorphic to the Seifert manifolds $\{-n; (o_1, 1) \}$ \cite[\S 7.2]{orlikSeifertManifolds1972}.
  \item $\tr(\gamma) = -2$: then there is a $\Z$-indexed family of conjugacy classes, given by $\{ - T^n : n \in \Z \}$. The corresponding mapping tori are diffeomorphic to $\{-n; (n_2, 2) \}$, in particular the mapping torus for $-I$ is $\{-2; (o_1, 0); (2, 1), (2, 1), (2, 1), (2, 1) \}$ \cite[\S 7.2]{orlikSeifertManifolds1972}.
  \item $|\tr(\gamma)| > 2$: then there is precisely one conjugacy class for each word of the form $R^{j_1}L^{k_1}\dots R^{j_I}L^{k_I}$ up to cyclic permutation of the sequence $(j_1, k_1, \dots, j_I, k_I)$ where $I \geq 1$ and $j_i, k_i$ are all positive integers. Here, $R = T$ and $L = T^{\top}$. These mapping classes are called hyperbolic, and they correspond to Anosov diffeomorphisms of the torus. By \cite[Lemma VI.31]{jacoLecturesThreemanifoldTopology1980}, since these do not have finite order and do not have trace $\pm 2$, their mapping tori are not Seifert manifolds. In fact, these manifolds are solvmanifolds, and any compact solvmanifold is finitely covered by some such mapping torus \cite[Thm. 5.3]{scottGeometries3Manifolds1983}.
\end{itemize}

We see that for $- 2 <  \tr(\gamma) < 2$ there are only 6 oriented diffeomorphism classes of mapping tori, while for $\tr(\gamma) = \pm 2$ there are $\Z$-indexed families of mapping tori. For $|\tr(\gamma)| > 2$ there is a description by cyclic sequences. As explained in Rmk. \ref{r-ccls-up-to-sign}, the $\SL_2$-skein modules of $M_{\pm \gamma}$ will have the same dimension. Table \ref{t-lowtrace-dims} gives the $\SL_2$-skein module dimensions for all mapping tori with $|\tr(\gamma)| \leq 2$, and Table \ref{t-length-2-dims} gives dimensions for mapping tori indexed by sequences of length 2.

Finally, we record for later use a basic computation about the homology of these manifolds.

\begin{lemma}
  \label{l-H_1-mapping-torus}
  The first homology of $M_\gamma$ is given by
  \[
    \H_1(M_\gamma) \cong \coker(\Id - \gamma) \oplus \Z
  \]
  where the above cokernel is for the induced map on homology $\Id - \gamma : \H_1(T^2) \to \H_1(T^2)$.
\end{lemma}

\begin{proof}
  From \cite[Example 2.48]{hatcherAlgebraicTopology2002}, there is a long exact sequence
  \[
    \dots \H_1(T^2) \xrightarrow{\Id - \gamma} \H_1(T^2)  \to \H_1(M_\gamma) \to H_0(T^2) \xrightarrow{\Id - \gamma} H_0(T^2) \dots
  \]
  and we note that $\Id - \gamma$ induces the zero map on $H_0(T^2)$. We can therefore extract the following short exact sequence
  \[
    0 \to \coker(\Id - \gamma) \to \H_1(M_\gamma) \to \Z \to 0
  \]
  and the result follows by the fact that $\Z$ is a projective $\Z$-module.
\end{proof}

%% file: 2-prelim/twisted-HH_0.tex
In this section we express the skein module of a mapping torus $M_\gamma$ as the Hochschild homology of the skein category of the 2-torus, twisted by the functor induced by $\gamma$. We begin by recalling skein categories of surfaces, then motivate and define the twisted Hochschild homology of skein categories, before stating an important result of \cite{gunninghamSkeinsTori} that allows us to give a direct sum decomposition of this twisted Hochschild homology.

The skein category was defined in  \cite{walkerTQFTs, johnson-freydHeisenbergpictureQuantumField2019}. We sketch the definition here.

\begin{defn}
  Let $\cA$ be a $k$-linear ribbon category, and $M$ a 3-manifold. For any oriented graph $\Gamma$ embedded in $M$, which admits half-edges meeting the boundary of $M$, we define a framing of $\Gamma$ to be a choice of a trivialization of the normal bundle to each (half) edge in the graph. A (half) edge in the graph equipped with a framing is called a \emph{ribbon}. An embedded oriented graph equipped with a framing and a cyclic ordering at each vertex will be called a \emph{ribbon tangle in $M$}. If a ribbon tangle does not meet the boundary of $M$ it will be called a \emph{ribbon graph}.

  An \emph{$\cA$-labelled ribbon tangle} is a ribbon tangle where every ribbon is labelled by an object of $\cA$, and vertices are labelled by morphisms from the ordered tensor product of the objects labelling ribbons oriented into the vertex to the ordered tensor product of those oriented outward.
\end{defn}

\begin{defn}
  The $k$-linear ribbon category $\mathrm{Rib}_\cA$ has
  \begin{itemize}
    \item Objects: finite unions of framed marked points in $I^2$ (that is, points together with a unit tangent vector) coloured by objects of $\cA$.
    \item Hom-spaces: $k$-vector spaces spanned by isotopy classes of ribbon tangles in $I^2 \times I$ which only meet the boundary at $I^2 \times \{0\}$ or $I^2 \times \{1\}$, with the framing of the ribbons agreeing with the framing vectors of points and the colourings agreeing. Composition is given by stacking in the third $I$ direction.
    \item Monoidal product given by disjoint union in the direction of the first interval copy of $I^2$.
    \item Duals given by reversing the direction of the framing of a point, with evaluation and coevaluation given by the cap and cup ribbons. The dual of a morphism (a ribbon tangle) is given by reversing the orientations of the edges, and dualizing the morphisms labelling the vertices.
    \item Braiding and twist given by the crossing and twisting of ribbons.
  \end{itemize}
  It is known, see for example \cite{turaevQuantumInvariantsKnots2010} for an overview, that there is a well-defined surjective and full functor $\mathrm{Rib}_\cA \to \cA$ of ribbon categories given by evaluating ribbon tangles, and the kernel of this functor defines what are called the \emph{$\cA$-skein relations}.
\end{defn}

\begin{defn}
  Let $\cA$ be a ribbon category. Then the \emph{$\cA$-skein module} $\Sk_\cA(M)$ of a 3-manifold $M$ is the $k$-vector space spanned by isotopy classes of ribbon graphs in $M$, modulo the $\cA$-skein relations which we apply in the interior of any embedded 3-ball $B^3 \cong I^3 \subseteq M$. An element of $\Sk_\cA(M)$ is called a \emph{skein}.
\end{defn}

In the case where $M = \Sigma \times I$ and the manifold has a distinguished interval direction, we can organize information about ribbon tangles modulo the skein relations (which we also call skeins) into the skein category.

\begin{defn}
  \label{d-skein-category}
  The $\cA$-\emph{skein category} of a surface $\Sigma$, denoted $\SkCat_\cA(\Sigma)$, has
  \begin{itemize}
    \item Objects: finite unions of framed marked points on $\Sigma$, coloured by objects of $\cA$. When a point is coloured with the monoidal unit in $\cA$ this is called the empty object, or \emph{empty skein}, and is denoted $\emptyset$.
    \item Hom-spaces: $k$-vector spaces spanned by isotopy classes of compatible $\cA$-labelled ribbon tangles in $\Sigma \times I$, modulo the $\cA$-skein relations (that is, the skein relations apply to sub-ribbon-graphs that lie in an embedded ball $B^3 \cong I^3 \subseteq \Sigma \times I$).
  \end{itemize}

  Given an object $x$, then $\End_{\SkCat_\cA(\Sigma)}(x)$ is an algebra with the product given by stacking. The \emph{$\cA$-skein algebra} of $\Sigma$ is the endomorphism algebra of $\emptyset \in \SkCat_\cA(\Sigma)$, denoted $\SkAlg_\cA(\Sigma)$.
\end{defn}

A typical choice of ribbon category is the category $\Rep^{\fd}_q(G)$. This is the category of finite-dimensional $k$-linear representations of the Drinfeld--Jimbo quantum group admitting a grading by the eigenspaces of the Cartan generators. While this Hopf algebra is defined over $\C(q)$, the explicit formula for its universal $R$-matrix obtained in \cite{Ros89AnalogueTheoremUniversal,KR90WeylGroupMultiplicative,LS90ApplicationsQuantumWeyl} shows that to consider the ribbon structure on its category of representations we must work over $k = \C(q^{\frac{1}{d}})$ where $d$ is the determinant of the symmetrized Cartan matrix of $G$. We refer readers to \cite[Chapters 8 - 10]{CP94GuideQuantumGroups} for details of the ribbon structure on $\Rep^{\fd}_q(G)$. In this case the skein module $\Sk_\cA(M)$ is denoted $\Sk_G(M)$. We denote by $\SkCat_G(\Sigma)$ the skein category and $\SkAlg_G(\Sigma)$ the skein algebra, for this choice of $\cA$.

\begin{notn}
  \label{n-framing}
  We will often denote skeins in $M_\gamma$ or in the skein category by giving unframed tangles in $T^2 \times I$, parameterized by $\{ (e^{2\pi i r}, e^{2\pi i s}, t): r, s, t \in [0, 1]\}$ where it is understood that $T^2 \times \{0\}$ and $T^2 \times \{1\}$ are identified via $\gamma$. We may simply give the tuple $(r, s, t) \in [0, 1]^3$ where the quotient is understood. The tangles we discuss will all admit, after a small isotopy, a projection to the cylinder $\{(e^{2 \pi i r}, e^{\pi i}, t) : r, t \in [0, 1]\}$ except perhaps at endpoints. The tangle edges can be canonically framed from this projection using the blackboard framing, and their endpoints framed in the direction of the first coordinate, and there is an obvious way to continue the framing of the projected edges to that of the endpoints. Issues of framing will not be important for our arguments but we mention this technical detail for completeness. The tangles and their endpoints are assumed to be coloured by the defining representation.
\end{notn}

Let us describe how to obtain the skein module of $M_\gamma$ from the twisted Hochschild homology of the skein category of $T^2$. Recall that we can associate to an algebra $E$ its Hochschild homology
\[
  \HH_0(E) = E/(ab - ba)
\]
which identifies the left and right actions of $E$ on itself. We can twist this when $\gamma$ is an automorphism of $E$.

\begin{defn}
  \label{d-twisted-HH_0}
  Let $E$ be an algebra and $\gamma$ an automorphism of $E$. Then the $(E, E)$-bimodule $E_\gamma$ has the underlying vector space $E$ with the left action given by left-multiplication
  \[
    a \rhd b = ab
  \]
  and the right action given by
  \[
    b \lhd a = b\gamma(a).
  \]
  We call $E_\gamma$ the \emph{$\gamma$-twisted bimodule}. The \emph{$\gamma$-twisted Hochschild homology} of $E$ is the space
  \[
    \HH_0^\gamma(E) := \HH_0(E, E_\gamma) = E/(ab - b\gamma(a))
  \]
  and the relations above are called the \emph{$\gamma$-commutator relations}.
\end{defn}

We regard categories as generalizing algebras, where algebras can be viewed as categories with a single object. Here, morphisms compose associatively and the appropriate generalization of Hochschild homology is given by coequalizing pre- and post-composition. This is given by the coend of the $\Hom$ bifunctor:
\[
  \HH_0(\cC) = \int^{x \in \cC} \Hom(x, x) = \left. \left( \bigoplus_{x \in \cC} \Hom(x, x) \right)\middle/ (f \circ g - g \circ f | g : x \to y, f : y \to x) \right. .
\]
Readers are referred to \cite[\S 1.2]{loregianCoEndCalculus2021} for the relevant background on coends. Notice that if $\cC$ has a single object $*$ then $\HH_0(\cC) = \HH_0(\End_\cC(*))$.

We consider the case when $\cC = \SkCat_{\cA}(\Sigma)$. It was proved in \cite[Lemma 5.5]{gunninghamFinitenessConjectureSkein2019} that for $\Sigma$ a closed surface,
\[
  \Sk_{\cA}(\Sigma \times S^1) \cong \HH_0(\SkCat_{\cA}(\Sigma))
\]
where the intuitive idea is that in identifying pre- and post-composition of morphisms in the skein category, the Hochschild homology construction gives a way to send skeins in $\Sigma \times I$ to their equivalence class on identifying the ends of the interval.

The argument of \cite{gunninghamFinitenessConjectureSkein2019} can be adapted to the situation where we twist by a mapping class $\gamma$, which acts on objects of the skein category, and also on morphisms by acting on $\Sigma \times I$ as $\gamma \times \Id$. This allows us to give a $\gamma$-twisted bifunctor $\Hom(\gamma(-), -)$. We think of skeins in $\Sigma \times_\gamma S^1$ as equivalence classes of skeins in $\Hom_{\SkCat_{\cA}(\Sigma)}(\gamma(x), x)$ where pre- and post-composition are identified, which gives a coend description of the skein module:

\begin{align}
  \Sk_{\cA}(\Sigma \times_\gamma I) & \cong \int^{x \in \SkCat_{\cA}(\Sigma)} \Hom(\gamma(x), x)\label{eq-coend-formula}                                                                                              \\
                                    & = \left. \left( \bigoplus_{x \in \SkCat_{\cA}(\Sigma)} \Hom(\gamma(x), x) \right) \middle/ (f \circ g - g \circ \gamma(f) | f : y \to x, g : \gamma(x) \to y) \right. \nonumber
  .
\end{align}

\begin{eg}
  It may be instructive to consider the examples illustrated in Figs. \ref{f-coend-example} and \ref{f-coend-example-twisted} for $\Sigma = T^2$.  Fig. \ref{f-coend-example} shows elements of different $\Hom$-spaces which should be identified on taking the coend of $\Hom(-, -)$, and Fig. \ref{f-coend-example-twisted} shows similar for the case twisted by $\gamma = S$.

  \begin{figure}
    \centering
    \includesvg[width=9cm]{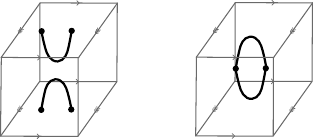}
    \caption{The figure depicts two morphisms in the skein category of $T^2$ which give the same skein on taking Hochschild homology. The left skein is the composition of the cup with the cap, an endomorphism of the two-points object; the right skein is the composition of the cap with the cup, an endomorphism of the empty object.}
    \label{f-coend-example}
  \end{figure}

  \begin{figure}
    \centering
    \includesvg[width=9cm]{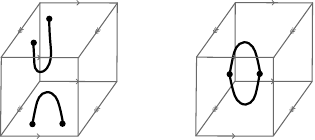}
    \caption{The figure depicts two morphisms in the skein category of $T^2$ which give the same skein on taking Hochschild homology twisted by $\gamma = S$. We depict $T^2$ as a quotient of $[0, 1]^2$, so that $S$ is rotation by a quarter turn, and let
      $x = \{(1/4, 1/4), (1/4, 3/4) \}$
      so that
      $\gamma(x) = \{(1/4, 1/4), (3/4, 1/4) \}$,
      and let $y = \gamma(y) = \emptyset$. Then, reading composition bottom to top, in the notation of (\ref{eq-coend-formula}), on the left we have a morphism $f \circ g : \gamma(x) \to x$ while on the right we have $g \circ \gamma(f) : y \to \gamma(y)$.}
    \label{f-coend-example-twisted}
  \end{figure}
\end{eg}

\begin{defn}
  For a $k$-linear category $\cC$, its \emph{free cocompletion} is the category $\Free(\cC) := [\cC^{\op}, \Vect]$, which is monoidal under Day convolution \cite{day1970construction}. We say that two categories $\cC, \cD$ are \emph{Morita equivalent} if and only if $\Free(\cC) \simeq \Free(\cD)$ (the usual notion of Morita equivalence for algebras is recovered for a one-object category).
\end{defn} 

The (twisted) Hochschild homology of $\cA$ is a Morita invariant, so in this paper we need only consider skein categories up to Morita equivalence and it suffices to replace $\SkCat_{\cA}(\Sigma)$ in (\ref{eq-coend-formula}) with $\Free(\SkCat_{\cA}(\Sigma))$. 

Recall that an object $x$ of a category is \emph{compact-projective} if $\Hom(x, -)$ preserves small colimits, and that the Yoneda embedding $\cJ : \cC \to \Free(\cC)$ exhibits $\cC$ as a full subcategory of compact-projective objects in $\Free(\cC)$. Recall also that the coend can be computed over any set of compact-projective generators. The following result gives a finite set of compact-projective generators for the free cocompletion of the $\GL_N$- and $\SL_N$-skein category of $T^2$.

\begin{thm}[{\cite[Thm. 1.5]{gunninghamSkeinsTori}}]
  \label{th-GJVY-decomp}
  Assume the parameter $q$ is generic. Then
  \begin{enumerate}
    \item The category $\Free(\SkCat_{\GL_N}(T^2))$ is generated by a single compact-projective object, given by the image under $\cJ$ of the empty skein, denoted by $\emptyset$.
    \item The category $\Free(\SkCat_{\SL_N}(T^2))$ admits $N$ compact-projective generators, given by image under $\cJ$ of the empty skein and the single point coloured by $V^{\ot n}, 0 < n < N$ (equivalently, $n$ points coloured by $V$) where $V$ is the defining representation of $U_q(\sl_N)$.
  \end{enumerate}
\end{thm}

Therefore, for $\Sigma = T^2, \cA = \Rep_q^{\fd}(G)$, we can take the coend (\ref{eq-coend-formula}) over the compact-projective generators of Thm. \ref{th-GJVY-decomp}. Since the generators are in the image of the Yoneda embedding, which is fully faithful, the $\Hom$-spaces appearing in (\ref{eq-coend-formula}) are $\Hom$-spaces in the skein category. Regarding each generator as a single point at the origin coloured by $V^{\ot n}$, we have that the generators are fixed by $\gamma$ and so the $\Hom$-spaces which appear are endomorphism algebras in the skein category.  

As described in \cite[\S 3]{Jor23LanglandsDualitySkein}, the $G$-skein category of $\Sigma$ is graded by $\H_1(\Sigma, Z(G)^{\vee})$, where $Z(G)^{\vee}$ is the Pontryagin dual of the centre of $G$. For $G = \SL_2$ this grading is known and was featured in the computations of \cite{carregaGeneratorsSkeinSpace2017,gilmerKauffmanBracketSkein2018,gilmerSkeinModuleProduct2019a,detcherryBasisKauffmanSkein2021}. Essentially, $Z(G)$ defines a grouplike central subalgebra of $U_q(\g)$, and using Schur's lemma this gives for each irreducible representation a central character, or element of $Z(G)^{\vee}$. This is compatible with the tensor structure of $\Rep_q^{\fd}(G)$ since the action of the centre is by grouplike elements. Because the skein relations take place locally in a ball the grading extends to a grading by $\H_1(\Sigma, Z(G)^{\vee})$ of the skein category. It follows that for $0 \leq n, m < N, n \neq m$ we have $\Hom_{\SkCat_{\SL_N}(\Sigma)}(V^{\ot m}, V^{\ot n}) = 0$, since $V^{\ot m}$ and $V^{\ot n}$ lie in different components of $\Rep_q^{\fd}(\SL_N)$ with respect to the grading by $\Z/N\Z$. 

Therefore, we have that the coend (\ref{eq-coend-formula}) splits as a direct sum of $\gamma$-twisted Hochschild homologies of endomorphism algebras:
\begin{equation}
  \label{eq-GL_N-decomp}
  \Sk_{\GL_N}(M_\gamma) \cong \HH_0^\gamma(\SkAlg_{\GL_N}(T^2))
\end{equation}
and
\begin{equation}
  \label{eq-SL_N-decomp}
  \Sk_{\SL_N}(M_\gamma) \cong \bigoplus_{n = 0}^{N-1} \HH_0^{\gamma}(\End_{\SkCat_{\SL_N}(T^2)}(V^{\ot n})).
\end{equation}

\begin{rmk}
  \label{r-SL_2-KBSM-equivalent}
  For a $k$-linear category, we denote by $\Kar(\cA) \subseteq \Free(\cA)$ the subcategory of all compact-projective objects of $\Free(\cA)$, also called the Karoubi completion of $\cA$. Note that $\Free(\Kar(\cA)) \simeq \Free(\cA)$. Then it follows from \cite[Thm. 3.27]{Coo23ExcisionSkeinCategories} that
  \[
    \Free(\SkCat_{\Kar(\cA)}(\Sigma)) \simeq \Free(\SkCat_{\cA}(\Sigma))
  \]
  so that $\SkCat_{\Kar(\cA)}(\Sigma)$ and $\SkCat_{\cA}(\Sigma)$ are Morita equivalent.

  Denoting by $\TL$ the full tensor subcategory of $\Rep_q^{\fd}(\SL_2)$ generated by the defining representation $V$ and its dual, it follows from complete reducibility (\cite{Ros90AnaloguesFormeKilling}, see \cite[Prop. 10.1.14]{CP94GuideQuantumGroups} for a textbook reference) that $\Kar(\TL) \simeq \Rep_q^{\fd}(\SL_2)$. Then for $G = \SL_2$, we may work with the $\TL$-skein category and the $\SL_2$-skein category interchangeably without changing our results.
  
  Quantum Schur-Weyl duality \cite{Jim86AnalogueMathfrakGl} says that $\End_{\TL}(V^{\ot n})$ is given by compositions and tensor products of the $R$-matrix for $V$ and its inverse, and this can be extended to general morphisms in $\TL$ (see e.g. \cite[Thm. 12.3.10]{CP94GuideQuantumGroups}). Then $\TL$ is equivalent to the category of oriented framed tangles modulo the Kauffman bracket skein relations, and since the defining representation is self-dual we can omit orientation data. This is called the Temperley-Lieb category in \cite{Tur90OperatorInvariantsTangles}, and we may call the associated skein categories and modules the Kauffman bracket skein category/module. We note that, since the skein module of a 3-manifold can be obtained from the Morita theory of the skein category of a Heegaard surface \cite[Cor. 1]{gunninghamFinitenessConjectureSkein2019}, this recovers the fact well-known to experts that the $\SL_2$ and Kauffman bracket skein modules are isomorphic for all 3-manifolds (and not just mapping tori).
\end{rmk}

%% file: 2-prelim/SkAlg.tex
In both of the decompositions (\ref{eq-GL_N-decomp}) and (\ref{eq-SL_N-decomp}) of the skein module, the twisted Hochschild homology of $\SkAlg_G(T^2)$ appears. In \S \ref{s-HH_0-SkAlg} we will give a decomposition of this space, which relies on a description of $\SkAlg_G(T^2)$ given by \cite{gunninghamQuantumCharacterTheory}. In this section, we recall this description and remark on the actions of $\SL_2(\Z)$ and the Weyl group on the skein algebra.

\begin{thm}[{\cite[Cor. 1.7]{gunninghamQuantumCharacterTheory}}]
  \label{th-GJVY_T^2-iso}
  There is an isomorphism
  \[
    \SkAlg_G(T^2) \cong (\cK_\omega[\Lambda \oplus \Lambda])^W.
  \]
\end{thm}

\begin{notn}
  \label{n-quantum-torus}
  In the above $\cK$ is the field of Notation \ref{n-fields} and $\Lambda$ is the weight lattice of $G$. There is the usual Cartan pairing $\la -, - \ra$ on $\Lambda$. We extend it to a nondegenerate skew pairing $\omega$ on $\Lambda \oplus \Lambda$ given by
\[
  \omega((r, s), (t, u)) = \la r, u \ra - \la s, t\ra.
\]
The algebra $\cK_\omega[\Lambda \oplus \Lambda]$ has as its underlying vector space the $\cK$-span of $\Lambda \oplus \Lambda$. If $x = (r, s) \in \Lambda \oplus \Lambda$ we denote the corresponding generator of the algebra by $m_{r, s}$ or $m_x$ and the the multiplication is given by
\[
  m_{x}m_{y} = q^{\frac{1}{2}\omega(x, y)}m_{x + y}.
\]
We will suppress much of this notation and denote by $A$ this $\cK$-algebra $\cK_\omega[\Lambda \oplus \Lambda]$. We denote by $W$ the Weyl group of $G$: it acts on $\Lambda \oplus \Lambda$ diagonally and this induces an action on $A$ by automorphisms.
\end{notn}

\begin{rmk}
  The $\SL_2$ case of Thm. \ref{th-GJVY_T^2-iso} was shown in \cite{frohmanSkeinModulesNoncommutative2000}. In \cite{frohmanSkeinModulesNoncommutative2000}, the basis $\{ m_{r, s} \}$ is named $\{ e_{r, s} \}$, the parameter $q^{1/2}$ is called $t$, and the algebra $\cK_{\omega}[\Lambda \oplus \Lambda]$ is called $\C_t[l, l^{-1}, m, m^{-1}]$ where $l = t^{-1}e_{1, 0}, m = t^{-1}e_{0, 1}$ so that $lm = t^2 ml$. Then the assignment $m_{r, s} \mapsto t^{-rs}l^rm^s, q^{1/2} \mapsto t$ establishes an algebra isomorphism $\cK_{\omega}[\Lambda \oplus \Lambda] \cong \C_t[l, l^{-1}, m, m^{-1}]$. The $\Z/2\Z$-actions on each algebra are equivalent, so that this restricts to an isomorphism of invariant subalgebras.
\end{rmk}

\begin{lemma}
  \label{l-invertible-iff-monomial}
  An element $a \in A$ is invertible if and only if $a = \alpha m_x$ for some $x \in \Lambda \oplus \Lambda, \alpha \in \cK^{\times}$.
\end{lemma}

\begin{proof}
  If $a = \alpha m_x$, then $a$ is invertible with inverse $\alpha^{-1} m_{-x}$. Conversely suppose that $a = \sum_{i=1}^n \alpha_i m_{x_i}$ is invertible with inverse $a^{-1} =  \sum_{i=1}^{n'} \beta_i m_{y_i}$, with nonzero coefficients $\alpha_i, \beta_i$. Then we can write $aa^{-1} = 1$ as
  \[
    \sum_{\substack{x_i, y_j\\x_i \neq - y_j}} \alpha_i \beta_j q^{\frac{1}{2}\omega(x_i, y_j)} m_{x_i + y_j} + (\sum_{\substack{x_i, y_j\\x_i = - y_j}} \alpha_i \beta_j - 1)m_0 = 0.
  \] 
  where in the $m_0$ term we used that $q^{\frac{1}{2}\omega(-y_i, y_j)} = 1$ since $\omega$ is skew. Linear independence of the $m_{x_i}$ implies that the coefficients of the first summation must vanish, but this is impossible, so the fist sum must be empty. This implies that $x_i = -y_j$ for all $i, j$, which is only possible if $n = n' = 1$ and so $a = \alpha_1 m_{x_1}$.
\end{proof}

\begin{lemma}
  \label{l-A-simple}
  The algebra $A$ is simple.
\end{lemma}

\begin{proof}
  Let $0 \neq I \subseteq A$ a two-sided ideal. By Lemma \ref{l-invertible-iff-monomial}, elements $\alpha m_x$ are invertible, so it suffices to exhibit such an element in $I$ to see that $1 \in I$ and $I = A$. 
  
  Let $a = \sum_{i = 1}^n \alpha_i m_{x_i} \in I$. We claim that $I$ contains an element of the form $\sum_{i=1}^{k} \beta m_{y_i}$ for all $1 \leq k \leq n$. Choosing one of the $m_{x_i}$ appearing in $a$, we have 
  \begin{align*}
    [m_{-x_i}, a] &= \sum_{1 \leq j \leq n} \alpha_j [m_{-x_i}, m_{x_j}]\\
     &= \sum_{1 \leq j \leq n} \alpha_j (q^{\frac{1}{2}\omega(-x_i, x_j)} - q^{\frac{1}{2}\omega(x_j, -x_i)}) m_{x_j - x_i}\\
      &= \sum_{\substack{1 \leq j \leq n\\ j \neq i}} \alpha_j (q^{\frac{1}{2}\omega(-x_i, x_j)} - q^{\frac{1}{2}\omega(x_j, -x_i)}) m_{x_j - x_i}
  \end{align*}
  where in the final equality we used that $\omega(-x_i, x_i) = \omega(x_i, -x_i) = 0$ since $\omega$ is skew. But since $I$ is a two-sided ideal, $[m_{x_i}, a] \in I$ and so we have shown our claim for $k = n-1$. Proceeding inductively, we can show that the claim holds for $k = 1$. Then $I$ contains an element of the form $\alpha m_x$, which is invertible, and so $I = A$.
\end{proof}

In the lattice $\Lambda \oplus \Lambda = \Lambda \ot \H_1(T^2)$, the two summands correspond to the two fundamental cycles in a torus, and $\SL_2(\Z)$ acts via its action on $\H_1(T^2)$: the mapping class $\gamma$ acts by $\Id \ot \gamma$. The diagonal action of $w \in W$ is by $w \ot \Id$ on $\Lambda \ot \H_1(T^2)$. Without further comment we will abuse notation and simply write $\gamma, w, \gamma w$ and so on for these automorphisms acting on the lattice $\Lambda \oplus \Lambda$ and hence on $A$.

\begin{lemma} 
  \label{l-pres-omega}
  The actions of $\SL_2(\Z)$ and $W$ on $\Lambda \oplus \Lambda$ both preserve the pairing $\omega$.
\end{lemma}

\begin{proof}
  Let $x = (r, s), y = (t, u) \in \Lambda \oplus \Lambda$. In case of $\SL_2(\Z)$, it suffices to show that the generators $S, T$ of $\SL_2(\Z)$ preserve $\omega$. We calculate
  \[
  \omega(S(r, s), S(t, u)) = \omega((-s, r), (-u, t)) = \la -s, t \ra - \la r, -u \ra = \omega((r, s), (t, u))
  \]
  and
  \[
    \omega(T(r, s), T(t, u)) = \omega((r + s, s), (t + u, u)) = \la r + s, u \ra - \la t + u, s\ra = \omega((r, s), (t, u))
  \]
  which establishes this. In the case of $w \in W$, we note that
  \[
    \omega(w(x), w(y)) = \omega(w(r, s), w(t, u)) = \la w(r), w(u) \ra - \la w(s), w(t) \ra = \la r, u \ra - \la s, t \ra = \omega(x, y)
  \]
  since $w$ preserves the Cartan pairing $\la-, -\ra$.
\end{proof}

It is easy to see that the actions of $\SL_2(\Z)$ and $W$ commute, so that $\SL_2(\Z)$ preserves the subspace of invariants $A^W$. From Lemma \ref{l-pres-omega}, we see that $\SL_2(\Z)$ acts by algebra automorphisms on $A$, and hence on $A^W$. It therefore makes sense to consider the bimodule $A^W_\gamma$ and the space $\HH_0^\gamma(A^W)$ as defined in Def. \ref{d-twisted-HH_0}. Finally, we make the following observation.

\begin{defn}
  An automorphism of an algebra is said to be \emph{outer} if it cannot be written as conjugation by some invertible element of the algebra.
\end{defn}

\begin{lemma}
  \label{l-W-outer}
  The action of $W$ on $A$ is by outer automorphisms.
\end{lemma}

\begin{proof}
  From Lemma \ref{l-invertible-iff-monomial}, the only invertible elements of $A$ are of the form $\alpha m_x$. Note that for $m_y \in A$, we have $(\alpha m_x) m_y (\alpha^{-1} m_{-x}) = m_y$, so conjugation by invertible elements preserves such elements. On the other hand, every nontrivial element of $W$ acts nontrivially on the set $\{ m_y : y \in \Lambda \oplus \Lambda \}$ since $W$ acts nontrivially on $\Lambda \oplus \Lambda$. So the action of nontrivial $w \in W$ on $A$ can never be written in terms of conjugation by an invertible element.
\end{proof}

%% file: 3-HH_0-SkAlg/outline.tex
In this section we consider the space $\HH_0^\gamma(\SkAlg_G(T^2))$, which is always a direct summand of $\Sk_G(M_\gamma)$. In \S \ref{s-HH_0_decomp} we give a direct sum decomposition of $\HH_0^\gamma(\SkAlg_G(T^2))$ into summands $\HH_0(A, A_{\gamma w})_{\Stab(w)}$ which are spaces of coinvariants for centralizers of the action of the Weyl group of $G$. In \S \ref{s-diff-cokernel} we give an expression for the spaces $\HH_0(A, A_{\gamma w})$ which makes it possible to deduce the dimensions of these spaces before taking coinvariants.

%% file: 3-HH_0-SkAlg/decomp.tex
In this section we give a direct sum decomposition for the space $\HH_0^\gamma(\SkAlg_G(T^2))$. We use the isomorphism of Thm. \ref{th-GJVY_T^2-iso} to write this as $\HH_0^\gamma(A^W)$. There is a Morita equivalence $\Phi$ of $A^W$ with the smash product $A \# W$, which is the $\cK$-algebra structure on $A \ot \cK W$ given on generators by
\[
  (a \# w_1) \cdot (b \# w_2) = aw_1(b) \# w_1 w_2
\]
and extended $\cK$-linearly, where $a \# w$ is notation for the tensor $a \ot w$ considered as an element in the smash product. The Morita equivalence $\Phi$ allows us to write the vector space we are interested in as
\[
  \HH_0^\gamma(\SkAlg_G(T^2)) \cong \HH_0(A^W, A^W_\gamma) \cong \HH_0(A\#W, \Phi(A^W_\gamma))
\]
and it is the rightmost expression which we will show admits a decomposition.

We begin by recalling the Morita equivalence and examining $\Phi(A^W_\gamma)$ in Prop. \ref{p-Morita-equiv-twisted-bimodule}, before using this to give the direct sum decomposition of $\HH_0^\gamma(\SkAlg_G(T^2))$, which is Thm. \ref{th-HH_0-decomp}. This allows us to refine the decompositions (\ref{eq-GL_N-decomp}, \ref{eq-SL_N-decomp}) in Cor. \ref{c-skmod-decomp}.

The Morita equivalence $\Phi$ follows from well-known arguments which we recall here for completeness. Throughout we assume that $|W|$ is invertible in $\cK$.

\begin{lemma}
  \label{l-Morita-std-lemma}
  Let $R$ be any ring, and $0 \neq e \in R$ an idempotent such that $ReR = R$. Then there is a Morita equivalence between $R$ and $S := eRe$ given by
  \begin{align*}
    \Phi : _{R}\Mod_R                     & \xleftrightarrow{\sim} _{S}\Mod_S : \Psi    \\
    M                                           & \mapsto eR \ot_R M \ot_R Re \\
    Re \ot_S N \ot_S eR & \mapsfrom N.
  \end{align*}
\end{lemma}

\begin{proof}
  This is a standard result. See e.g. \cite[Ex. 21.6]{AF74RingsCategoriesModules} for the equivalence ${}_R \Mod \simeq {}_S \Mod$ of categories of left modules implemented by $eR$ and $Re$. It is known that $_{R} \Mod \simeq _{S}\Mod$ if and only if $\Mod_{R} \simeq \Mod_{S}$. Applying both equivalences simultaneously gives the claimed equivalence for the category of bimodules. 
\end{proof}

\begin{lemma}
  \label{l-smash-product-simple}
  Let $R$ be a ring with unit, and $G$ a group acting on $R$ by outer automorphisms. Then if $R$ is simple it follows that $R \# G$ is simple.
\end{lemma}

\begin{proof}
  This is given in \cite[Thm. 2.3]{Mon80FixedRingsFinite} in the case of the skew product $R * G$, but the proof also works in the case of the smash product.
\end{proof}

\begin{prop}
  There is a Morita equivalence
  \begin{align*}
    \Phi : _{A^W}\Mod_{A^W}                     & \xleftrightarrow{\sim} _{A\#W}\Mod_{A\#W} : \Psi    \\
    M                                           & \mapsto {}_{A\# W}A \ot_{A^W} M \ot_{A^W} A^W_{A\# W} \\
    {}_{A^W}A^W \ot_{A\# W} N \ot_{A\# W} A_{A^W} & \mapsfrom N
  \end{align*}
  where ${}_{A^W} A^W_{A\# W}$ is the bimodule on $A^W$ with the obvious left action and the right action
  \[
    a \lhd (b \# w) = ae(b)
  \]
  and ${}_{A\# W}A \ot_{A^W}$ is the bimodule on $A$ with the right action given by the action of a subring and the left action given by
  \[
    (b \# w) \rhd a = bw(a).
  \]
\end{prop}

\begin{proof}
  By Lemmas \ref{l-A-simple}, \ref{l-W-outer}, $A$ is simple with $W$ acting by outer automorphisms. Then by Lemma \ref{l-smash-product-simple}, $A \# W$ is simple, and we can apply Lemma \ref{l-Morita-std-lemma} using any idempotent of $A \# W$. 
  
  Let $f = \sum_{w \in W} w$ and $e = \frac{1}{|W|} f$. Then $e$ is an idempotent in $\cK W$ and $1 \# e$ is an idempotent in $A \# W$. Moreover, it is easily checked that
  \begin{equation*}
    (1 \# e) (A \# W) = (1 \# e) (A \# W) (1 \# e) = A^W \# e \cong A^W
  \end{equation*}
  where the last isomorphism is an isomorphism of $(A^W, A\#W)$-bimodules with the bimodule structures given in the statement of the proposition. Similarly,
  \begin{equation*}
    (A \# W) (1 \# e) = A \# e \cong A
  \end{equation*}
  as an $(A \# W, A^W)$-bimodule under the stated actions.
\end{proof}


\begin{prop}
  \label{p-Morita-equiv-twisted-bimodule}
  We have $\Phi(A^W_\gamma) \cong A_\gamma \# W$ where $A_\gamma \# W$ is the bimodule on $A \ot \cK W$ with actions
  \begin{align*}
    b \# w_2 \rhd a \ot w_1 & = bw_2(a) \ot w_2w_1          \\
    a \ot w_1 \lhd b \# w_2 & = a\gamma(w_1(b)) \ot w_1w_2.
  \end{align*}
\end{prop}

\begin{proof}
  It is straightforward to check that this is indeed an $A \# W$-bimodule structure. Then it suffices to show that $\Psi(A_\gamma \# W) \cong A^W_\gamma$ under the equivalence. 

  We first observe that for any right $A \# W$-module $M$, with the right action denoted $\lhd_1$, we have an isomorphism of right $A \# W$-modules
  \begin{equation}
    \label{eq-Morita-equiv-1}
    M \ot_{A \# W} A_\gamma \# W \cong M
  \end{equation}
  where the $A \# W$-module structure on the right hand side, denoted $\lhd_2$, is given by $m \lhd_2 (a \# w) = m \lhd_1 (\gamma(a) \# w)$.

As underlying vector spaces, we have
  \begin{align}
    A^W \ot_{A \# W} A_\gamma \# W \ot_{A \# W} A & \cong A^W \ot_{A \# W} A \label{Morita-congruences-1} \\
                                                & \cong A^W \ot_{A \# W} A \# W \ot_{\cK W} \cK \label{Morita-congruences-2}                                            \\
                                                & \cong A^W \ot_{\cK W} \cK \label{Morita-congruences-3}                                                                         \\
                                                & \cong A^W. \label{Morita-congruences-4}
  \end{align}
  Line (\ref{Morita-congruences-1}) is an application of (\ref{eq-Morita-equiv-1}). Line (\ref{Morita-congruences-2}) applies the obvious isomorphism $A \cong A \# W \ot_{\cK W} \cK$. Line (\ref{Morita-congruences-3}) is an untwisted version of (\ref{eq-Morita-equiv-1}), and line (\ref{Morita-congruences-4}) follows since the $W$-action on $A^W$ is trivial.

  Now reversing the isomorphisms (\ref{Morita-congruences-1} - \ref{Morita-congruences-4}), we have that $a \in A^W$ is sent to $a \ot 1 \ot 1 \ot 1$ where we abuse notation and choose representatives of equivalence classes. Then the right action of $A^W$ is given by $(a \ot 1 \ot 1 \ot 1) \lhd b = a \ot 1 \ot 1 \ot b$ which maps under the equivalences (\ref{Morita-congruences-1}, \ref{Morita-congruences-2}) to $a \ot b \ot 1 \ot 1$ and then, under (\ref{Morita-congruences-3}, \ref{Morita-congruences-4}) this maps to $a\gamma(b) \in A^W$ using the right action of $A \# W$ on $A^W$ that arises from line (\ref{Morita-congruences-1}). So the right $A^W$-module structure on $A^W$ is twisted by $\gamma$.

  It is clear that the left action is unaffected, and so it follows that $A_{\gamma} \# W$ is sent to $A^W_\gamma$ under the Morita equivalence $\Psi$.
\end{proof}

Having understood the image of $A^W_\gamma$ under $\Phi$, we are now ready to give a decomposition of the Hochschild homology.

\begin{thm}
  \label{th-HH_0-decomp}
  The canonical map
  \begin{align*}
    \HH_0(A \# W, A_{\gamma} \# W) & \to \bigoplus_{[w] \in \conj{W}} \HH_0(A, A_{\gamma w})_{\Stab(w)} \\
    [a \# w]                       & \mapsto [a] \in \HH_0(A, A_{\gamma w})_{\Stab(w)}
  \end{align*}
  is an isomorphism of vector spaces.
\end{thm}

\begin{proof}
  A general commutator relation in $\HH_0(A \# W, A_\gamma \# W)$ is of the form
  \[
    [a \# w_1, b \ot w_2] = aw_1(b) \ot w_1 w_2 - b\gamma(w_2(a)) \ot w_2 w_1
  \]
  where we understand the right argument of the bracket to be in the bimodule, and the left argument to be in the algebra acting. We are interested in the quotient $A_\gamma \# W/\la S \ra$ where
  \[
    S = \{ aw_1(b) \ot w_1 w_2 - b\gamma(w_2(a)) \ot w_2 w_1 : w_1, w_2 \in W; a, b \in A\}.
  \]

  We claim that $\la S \ra = \la S_W \cup S_A \ra$ where $S_W = \{[a \# g^{-1}w, 1 \ot g] : w, g \in W; a \in A\}$ and $S_A = \{[a \# 1, b \ot w]: w \in W; a, b \in A \}$. The inclusion $S \supseteq S_W \cup S_A$ is clear. For the reverse inclusion of linear spans, take any general relation
  \[
    aw_1(b) \ot w_1 w_2 - b\gamma(w_2(a)) \ot w_2 w_1.
  \]
  Now, let $w_2' = w_1 w_2$ so that $w_2 = w_1^{-1}w_2'$. Then the relation has the form
  \[
    aw_1(b) \ot w_2' - b\gamma(w_2(a)) \ot w_1^{-1}w_2' w_1.
  \]
  Now let $a' = w_1^{-1}(a)$ so that $a = w_1(a')$. Then the relation has the form
  \[
    w_1(a'b) \ot w_2' - b\gamma(w_1^{-1}w_2'w_1(a')) \ot w_1^{-1}w_2' w_1.
  \]
  Considering
  \[
    R_1 = [a' \# 1, b \ot w_1^{-1}w_2'w_1] = a'b \ot w_1^{-1}w_2'w_1 - b\gamma(w_1^{-1}w_2'w_1(a')) \ot w_1^{-1}w_2'w_1 \in S_A
  \]
  and
  \[
    R_2 = [a'b \# w_1^{-1}w_2', 1 \ot w_1] = a'b \ot w_1^{-1}w_2'w_1 - \gamma(w_1(a'b)) \ot w_2' \in S_W
  \]
  and
  \[
    R_3 = [w_1(a'b) \# w_2', 1 \ot 1] = w_1(a'b) \ot w_2' - \gamma(w_1(a'b)) \ot w_2' \in S_W
  \]
  then we see that our original relation is $R_1 - R_2 + R_3$
  so lies in the desired span.

  Taking the quotient of $A_\gamma \# W$ first by the span of $S_A$ yields
  \[
    A_\gamma \# W / \la S_A \ra \cong \bigoplus_{w \in W} \HH_0(A, A_{\gamma w})
  \]
  since the relations in $S_A$ fix the $W$-coordinate and on the $A$-coordinate are the $\gamma w$-twisted Hochschild relations. Further taking the quotient by the action of $S_W$, we see that
  \[
    [a \# g^{-1}w, 1 \ot g] = a \ot g^{-1}wg - \gamma(g(a))\ot w
  \]
  so that each such relation identifies $\HH_0(A, A_{\gamma w})$ and $\HH_0(A, A_{\gamma g^{-1}wg})$, twisted by the automorphism given by $\gamma g$. So it suffices to consider the space
  \[
    \bigoplus_{[w] \in \conj{W}} \HH_0(A, A_{\gamma w})
  \]
  where the sum ranges over conjugacy classes.

  There are further relations in each summand above. Suppose we have chosen a representative $w$ of the conjugacy class $[w]$, and there is $w'$ with $g, k$ such that both $g, k$ conjugate $w'$ to $w$. This will be the case if and only if $g = sk$ for some $s \in \Stab(w)$. Then the remaining relations in $\HH_0(A, A_{\gamma w})$ are those from $S_W$ where $g \in \Stab(w)$, in which case the relations $[a \# g^{-1}w, 1 \ot g]$ are
  \[
    \{a - \gamma(g(a)) : g \in \Stab(w)\}.
  \]
  We call these the $\gamma$-centralizer relations. These resemble the relations for taking coinvariants for the action of the centralizer:
  \[
    \{a - g(a) : g \in \Stab(w)\}
  \]
  but where the group action is twisted by $\gamma$.

  So far we have established an isomorphism
  \[
    A_\gamma \# W/\la S_A \cup S_W \ra \cong \bigoplus_{[w] \in \conj{W}} \HH_0(A, A_{\gamma w})/\la a - \gamma(g(a)) : g \in \Stab(w) \ra
  \]
  and it remains to show that the $\gamma$-centralizer relations are equivalent to the relations of taking coinvariants for the centralizer. We see that the $\gamma$-centralizer relations imply that for $s, t \in \Stab(w)$ we have
  \[
    \gamma(s(a)) \sim a \sim \gamma(t(a))
  \]
  and then since $e \in \Stab(w)$ and $\gamma$ is an automorphism, this implies $a \sim s(a)$ by taking $t = e$ and inverting $\gamma$. Conversely, suppose we take the quotient of $\HH_0(A, A_{\gamma w})$ by the relations $a \sim s(a)$ for $s \in \Stab(w)$. Then we have that
  \[
    w(a) \sim a \sim s(a)
  \]
  and applying $\gamma$,
  \[
    \gamma(w(a)) \sim \gamma(s(a))
  \]
  but we already have $a \sim \gamma(w(a))$ in $\HH_0(A, A_{\gamma w})$. So $a \sim \gamma(s(a))$ and we recover the $\gamma$-centralizer relations. It follows that taking the quotient of $\HH_0(A, A_{\gamma w})$ by the $\gamma$-centralizer relations is equivalent to taking the quotient by taking coinvariants for the usual action of the centralizer. This concludes the proof.
\end{proof}

\begin{cor}
  \label{c-skmod-decomp}
  For $G = \GL_N$, we have a decomposition of the skein module of $M_\gamma$ as
  \[
    \Sk_{\GL_N}(M_\gamma) \cong \bigoplus_{[w] \in \conj{W}} \HH_0(A, A_{\gamma w})_{\Stab(w)}.
  \]
  For $G = \SL_N$, we have the decomposition
  \[
    \Sk_{\SL_N}(M_\gamma) \cong \bigoplus_{n = 1}^{N-1} \HH_0^{\gamma}(\End_{\SkCat_{\SL_N}(T^2)}(V^{\ot n})) \oplus \bigoplus_{[w] \in \conj{W}} \HH_0(A, A_{\gamma w})_{\Stab(w)}.
  \]
\end{cor}

\begin{rmk}
  \label{r-ccls-up-to-sign}
  In this setting, it is clear that when $W = \Z/2\Z$ acts by $\pm 1$, the skein modules for $\gamma$ and $-\gamma$ will have the same dimensions. This is specific to the $\SL_2$-skein module, and for different groups we may not make this simplification.
\end{rmk}

%% file: 3-HH_0-SkAlg/diff-cokernel.tex
In this section we will express $\HH_0(A, A_{\gamma w})$ in terms of the cokernel of $\Id - \gamma w$. This gives a description of $\dim \HH_0(A, A_{\gamma w})$ in terms of the invariant factors of $\Id - \gamma w$. We illustrate how this fits into the skein module computation in the simple example of $G = \GL_1$. First let us introduce some notation.

\begin{notn}
  We will abbreviate $\Lambda \oplus \Lambda$ as $\Lambda^2$. There is a map of $\Z$-modules $\Id - \gamma w : \Lambda^2 \to \Lambda^2$, which extends to an endomorphism of $A$ as a $\cK$-vector space which we denote identically. As well as the $\cK$-vector space $A$, we can also consider the $\Q$-vector space $\Lambda^2_\Q = \Lambda^2 \ot_\Z \Q$ which is the linearization of $\Lambda^2$. Then $\Id - \gamma w$ extends to a $\Q$-linear map $(\Id - \gamma w)_\Q$, and we denote the kernel by $K_\Q$. We refer to the orthogonal complement with respect to $\omega$ as $K^{\perp}_\Q$. The restriction back to the lattice is denoted $K = K_\Q \cap \Lambda^2, K^\perp = K^{\perp}_\Q \cap \Lambda^2$. Note that these objects have been defined such that $\dim_\Q V_\Q = \rk_\Z V$, for $V$ a submodule of $\Lambda^2$.
\end{notn}

\begin{lemma}
  \label{l-K-perp}
  Suppose that $m_x \in A$ does not lie in $K^\perp$. Then $m_x$ is zero under the $\gamma w$-commutator relations.
\end{lemma}

\begin{proof}
  Suppose $m_x \notin K^\perp$. Then there exists $m_y \in K$ such that $\omega(x, y) \neq 0$, and notice that $m_y = \gamma w (m_y)$ since $m_y$ is in the kernel of $\Id - \gamma w$. Then consider the relation
  \[
    [m_y, m_{-y}m_x]_{\gamma w} = (1 - q^{\omega(x, y)})m_x
  \]
  which is easily verified. Then, in the quotient by the twisted commutator relations, $m_x$ must vanish.
\end{proof}

Therefore we have that $\HH_0(A, A_{\gamma w}) \cong \cK \la K^{\perp} \ra /(\gamma w\text{-commutators})$.

\begin{rmk}
  \label{r-when-K^perp-zero}
  When $K = 0$, then $K^\perp = \Lambda^2$. This is the case when $\det(\Id - \gamma w) \neq 0$. Notice that $\det(\Id - \gamma w)$ is the characteristic polynomial of $\gamma w$ evaluated at 1. In the case of a 2 by 2 matrix with determinant 1, this is $2 - \tr(\gamma w)$. So we will only need to consider $K^\perp$ a proper sublattice when $\tr(\gamma w) = 2$, that is, when $\gamma w$ corresponds to a conjugacy class of a matrix $T^n$ for some $n$.
\end{rmk}

We now give a useful change of basis in $K^\perp$. To describe the change of basis, we will consider the linearization $\Lambda^2_\Q$ etc. From Lemma \ref{l-pres-omega}, we see that $\gamma w$ preserves $\omega$, i.e. $\omega(\gamma w(x), \gamma w(y)) = \omega(x, y)$.

Then let $ x = y - \gamma w(y) \in \Img(\Id - \gamma w)_\Q$ and $k \in K_\Q$, so $k = \gamma w(k)$. We have
\[
  \omega(x, k) = \omega(y - \gamma w(y), k) = \omega(y, k) - \omega(\gamma w(y), k) = \omega(y, k) - \omega(\gamma w(y), \gamma w(k)) = \omega(y, k) - \omega(y, k) = 0
\]
so we see that $\Img(\Id - \gamma w)_\Q \subseteq K_\Q^\perp$. Moreover, we have that
\[
  \dim_\Q \Img(\Id - \gamma w)_\Q = \dim_\Q\Lambda^2_\Q - \dim_\Q K_\Q = \dim_\Q K_\Q^\perp
\]
where the first equality is the rank-nullity theorem and the second is from the fact that a subspace and its orthogonal complement split a vector space as a direct sum. Then it follows that 
\begin{equation}
  \label{eq-K-perp-image}
  \Img(\Id - \gamma w)_\Q = K_\Q^\perp.
\end{equation}
Notice that (\ref{eq-K-perp-image}) would not necessarily hold if we had not considered linearizations: consider $\gamma = T^{-2}$ for example, to see that $\Img(\Id - \gamma)$ may be a proper sublattice of $K^\perp$.

We see from (\ref{eq-K-perp-image}) that $\Lambda^2_\Q$ surjects onto $K_\Q^\perp = \Img(\Id - \gamma w)_\Q \cong \Lambda^2_\Q/K_\Q$. Given $x \in K_\Q^\perp$, its preimage is a $K_\Q$-coset. Let $(\Id - \gamma w)^{-1}x$ be a choice of representative, and notice that for $k \in K_\Q$ we have
\[
  \omega((\Id - \gamma w)^{-1}x + k, x) = \omega((\Id - \gamma w)^{-1}x, x) + \omega(k, x) = \omega((\Id - \gamma w)^{-1}x, x)
\]
so that the quantity $\omega((\Id - \gamma w)^{-1}x, x)$ is well-defined for $x \in K_\Q^{\perp}$, in particular for $x \in K^\perp \subseteq K_\Q^\perp$.

\begin{prop}
  \label{p-K-perp-renormalize}
  The change of basis
  \[
    m_x \mapsto \tilde{m}_x = q^{\frac{1}{2}\omega((\Id - \gamma w)^{-1}x, x)}m_{x}
  \]
  of $K^{\perp}$ is such that the commutator relations are all proportional to
  \[
    \tilde{m}_{x + y} - \tilde{m}_{x + \gamma w(y)}.
  \]
\end{prop}

\begin{proof}
  Consider the general commutator, which has the form
  \begin{equation}
    \label{eq-generic-twisted-commutator}
    [m_x, m_y]_{\gamma w} = q^{\frac{1}{2}\omega(x, y)}m_{x + y} - q^{\frac{1}{2}\omega(y, \gamma w(x))}m_{y + \gamma w(x)}.
  \end{equation}
  A renormalization $\tilde{m}_z = q^{-f(z)}m_z$ will imply the commutators have the desired form if the coefficients appearing in (\ref{eq-generic-twisted-commutator}) after this substitution have the same absolute value and differ just by a sign. This is equivalent to saying that
  \begin{equation}
    \label{eq-renormalization-constraint}
    f(x + y) - f(\gamma w(x) + y) = \frac{1}{2}\omega(y, \gamma w(x)) - \frac{1}{2}\omega(x, y).
  \end{equation}

  We claim that $f(z) = -\frac{1}{2}\omega((\Id - \gamma w)^{-1}z, z)$ satisfies (\ref{eq-renormalization-constraint}).

  To prove this, we will use repeatedly the fact that $\gamma w$ preserves the form $\omega$, and moreover the identity
  \begin{equation}
    \label{eq-1-P-inv-identity}
    (\Id - \gamma w)^{-1}\gamma w = \gamma w(\Id - \gamma w)^{-1} = (\Id - \gamma w)^{-1} - \Id.
  \end{equation}
  (We note that, where we interpret $(\Id - \gamma w)^{-1}$ as involving a choice of preimage, then (\ref{eq-1-P-inv-identity}) only holds up to a $K_\Q$-coset. The identity will only ever be applied to a single side of $\omega$ against either $x$ or $\gamma w(x)$ for $x \in K_\Q^\perp$. Then we note that, $k \in K_\Q$ means $\gamma w(k) = k$ and so $\omega(\gamma w(x), k) = \omega(\gamma w(x), \gamma w(k)) = \omega(x, k) = 0$ so that $\gamma w(x) \in K_\Q^{\perp}$. Then, by a similar remark to above, we see that it will be valid to apply (\ref{eq-1-P-inv-identity}) in our circumstances.)

  Let us expand the left hand side of (\ref{eq-renormalization-constraint}).

  \begin{align*}
      f(x &+ y) - f(\gamma w(x) + y)\\
      = &-\frac{1}{2}\omega((\Id - \gamma w)^{-1}(x + y), (x + y)) + \frac{1}{2}\omega((\Id - \gamma w)^{-1}(\gamma w(x) + y), \gamma w(x) + y)\\
      = &-\frac{1}{2}\omega((\Id - \gamma w)^{-1}x, x) -\frac{1}{2}\omega((\Id - \gamma w)^{-1}x, y) -\frac{1}{2}\omega((\Id - \gamma w)^{-1}y, x)\\
      &-\frac{1}{2}\omega((\Id - \gamma w)^{-1}y, y) + \frac{1}{2}\omega((\Id - \gamma w)^{-1}\gamma w(x) ,\gamma w(x)) + \frac{1}{2}\omega((\Id - \gamma w)^{-1} \gamma w(x),y)\\
      &+ \frac{1}{2}\omega((\Id - \gamma w)^{-1} y, \gamma w(x)) + \frac{1}{2}\omega((\Id - \gamma w)^{-1} y, y)\\
      = &-\frac{1}{2}\omega((\Id - \gamma w)^{-1}x, y) -\frac{1}{2}\omega((\Id - \gamma w)^{-1}y, x) + \frac{1}{2}\omega((\Id - \gamma w)^{-1} \gamma w(x),y)\\
      &+ \frac{1}{2}\omega((\Id - \gamma w)^{-1} y, \gamma w(x))
  \end{align*}
  where in the final equality we used that 
  \[
    \frac{1}{2}\omega((\Id - \gamma w)^{-1}\gamma w(x) ,\gamma w(x)) = \frac{1}{2}\omega(\gamma w(\Id - \gamma w)^{-1}x ,\gamma w(x)) = \frac{1}{2}\omega((\Id - \gamma w)^{-1}x ,x)
  \]
  from (\ref{eq-1-P-inv-identity})
  and that $\gamma w$ preserves $\omega$. Now we can continue, applying the identity (\ref{eq-1-P-inv-identity}) again in the third summand.
  \begin{align*}
      f(x &+ y) - f(\gamma w(x) + y)\\
      = &-\frac{1}{2}\omega((\Id - \gamma w)^{-1}x, y) -\frac{1}{2}\omega((\Id - \gamma w)^{-1}y, x) + \frac{1}{2}\omega((\Id - \gamma w)^{-1} x,y) - \frac{1}{2}\omega(x,y)\\
      & + \frac{1}{2}\omega((\Id - \gamma w)^{-1} y, \gamma w(x))\\
      = &-\frac{1}{2}\omega((\Id - \gamma w)^{-1}y, x) - \frac{1}{2}\omega(x,y) + \frac{1}{2}\omega((\Id - \gamma w)^{-1} y, \gamma w(x))\\
      = &-\frac{1}{2}\omega(\gamma w(\Id - \gamma w)^{-1}y, \gamma w(x)) - \frac{1}{2}\omega(x,y) + \frac{1}{2}\omega((\Id - \gamma w)^{-1} y, \gamma w(x))\\
      = &-\frac{1}{2}\omega((\Id - \gamma w)^{-1}y, \gamma w(x)) +\frac{1}{2}\omega(y, \gamma w(x)) - \frac{1}{2}\omega(x,y) + \frac{1}{2}\omega((\Id - \gamma w)^{-1} y, \gamma w(x))\\
      = &\hspace{12.5pt}\frac{1}{2}\omega(y, \gamma w(x)) - \frac{1}{2}\omega(x,y)
  \end{align*}
  using that $\gamma w$ preserves $\omega$ in the third equality, and identity (\ref{eq-1-P-inv-identity}) in the fourth. This establishes the required identity.
\end{proof}

We can use Prop. \ref{p-K-perp-renormalize} to obtain a description of the space $\HH_0(A, A_{\gamma w})$.

\begin{cor}
  \label{c-dim-HH_0-before-stab}
  The linear map
  \begin{align*}
    \HH_0(A, A_{\gamma w}) & \to \cK[\coker(\Id - \gamma w)_{\mathrm{tors}}]\\
    [m_x]                  & \mapsto [x]
  \end{align*}
  is an isomorphism of vector spaces, where the right hand side denotes the vector space supported on the torsion subgroup of the cokernel of $\Id - \gamma w : \Lambda^2 \to \Lambda^2$. Let $a_i^w$ be the invariant factors of $\Id - \gamma w$, and let $r_w$ be the rank of this map. The dimension of this space is
  \[
    \dim \HH_0(A, A_{\gamma w}) = \prod_{i=1}^{r_w} a_i^w.
  \]
\end{cor}

\begin{proof}
  Notice that under Prop. \ref{p-K-perp-renormalize}, the $\gamma w$-commutator relations are equivalent (re-indexing) to $\{\tilde{m}_x - \tilde{m}_{x + (\Id - \gamma w)(y)}\}.$ But this is the same as taking the quotient of $K^\perp$ (as a $\Z$-module) by the submodule $\Img(\Id - \gamma w)$.

  We consider $\Id - \gamma w : \Lambda^2 \to \Lambda^2$ as a map of $\Z$-modules. Then since $\Lambda^2$ is finitely generated, the cokernel of $\Id - \gamma w$ is isomorphic to $\Z^l \oplus \bigoplus_i \Z / a_i$ where $a_i$ are the invariant factors of the map. The number of torsion summands is the rank of the map, and $l$ is the corank of the map.

  We recall that $\Img(\Id - \gamma w) \subseteq K^\perp$ and, since by (\ref{eq-K-perp-image}) we have $\Img(\Id - \gamma w)_\Q \cong K^\perp_\Q$, we then have that $\rk_\Z \Img(\Id - \gamma w) = \dim_\Q \Img(\Id - \gamma w)_\Q = \dim_\Q K^\perp_\Q = \rk_\Z K^\perp$. Then it follows that the torsion part of $\coker(\Id - \gamma w : \Lambda^2 \to \Lambda^2)$ is precisely $\coker(\Id - \gamma w : \Lambda^2 \to K^\perp)$. This establishes the first statement. The cardinality of the abelian group is given by taking products, and it is this cardinality which gives the dimension of the vector space over $\cK$, giving the second statement.
\end{proof}

From Cor. \ref{c-dim-HH_0-before-stab}, using the decomposition of Cor. \ref{c-skmod-decomp}, then to understand $\dim \HH_0^\gamma(\SkAlg_G(T^2))$ it suffices to understand the centralizers of the Weyl group and their orbits. We finish this section by considering the simple case of $G=\GL_1$, where the Weyl group is trivial.

We recall that it was shown in \cite{Prz98AnalogueFirstHomology} that there is an isomorphism $\Sk_{\GL_1}(M) \cong \C(q)[\H_1(M)_\mathrm{tors}]$ for $q$ generic. In the case of a mapping torus $M_\gamma$, we mentioned in Lemma \ref{l-H_1-mapping-torus} how to compute $\H_1(M_\gamma)$. To compute the torsion subgroup, it suffices to find the torsion of
\[
  \coker(\Id - \gamma) = \Z \oplus \Z /\Img(\Id - \gamma).
\]
When $\Id - \gamma$ is invertible, then the cardinality of its image is measured by its determinant. We have that $\Id - \gamma$ is invertible when $\det(\Id - \gamma)$ is nonzero, but since $\det(\Id - \gamma)$ is the characteristic polynomial of the 2 by 2 matrix $\gamma$ evaluated at 1, and $\gamma$ has determinant 1, it follows that $\det(\Id - \gamma) = 2 - \tr(\gamma)$ and so $\Id - \gamma$ is invertible for all matrices except the shears. For the shears $T^n$ we have that $\Img(\Id - T^n) = \Z / n\Z \oplus \Z$. It follows that
\begin{equation}
  \label{eq-GL_1-dims}
  \dim \Sk_{\GL_1}(M_\gamma) = |\H_1(M_\gamma)_\mathrm{tors}| = \begin{cases}
    1                                        & \gamma = \Id           \\
    |n|                                      & \gamma = T^n, n \neq 0 \\
    |\det(\Id - \gamma)| = |\tr(\gamma) - 2| & \text{otherwise.}
  \end{cases}
\end{equation}
We note that this result can be recovered independently of \cite{Prz98AnalogueFirstHomology} by our methods.

\begin{thm}
  \label{th-GL_1-dims}
  The dimension of the $\GL_1$-skein module of $M_\gamma$ is given as in (\ref{eq-GL_1-dims}).
\end{thm}

\begin{proof}
  Here the Weyl group is trivial: the sum in Cor. \ref{c-skmod-decomp} consists of a single summand $\HH_0(A, A_\gamma)$ and all that remains is to take the quotient of $A$ by the $\gamma$-commutator relations.

  By lemma \ref{l-K-perp} we need only consider the quotient of $K^\perp$, which only differs from $A$ in the case of the shear matrices, by Rmk. \ref{r-when-K^perp-zero}. In the shear case $\gamma = \Id$ or $\gamma = T^n$, the dimension follows from the description of Cor. \ref{c-dim-HH_0-before-stab}.

  In the non-shear cases we want to consider $A$ modulo the image of $\Id - \gamma$. But in this case, since $\Id - \gamma$ has nonzero determinant, then the cardinality of this cokernel is given by $|\det(\Id - \gamma)|$.
\end{proof}

%% file: 4-SL_2-dims/outline.tex
In this section we focus on the specific case of $G = \SL_2$. We recall the decomposition (\ref{eq-SL_N-decomp}):
\[
  \Sk_{\SL_2}(M_\gamma) = \HH_0^\gamma(\End_{\SkCat_{\SL_2}(T^2)}(V)) \oplus \HH_0^\gamma(\SkAlg_{\SL_2}(T^2)).
\]
Here, $V$ denotes the object of the skein category consisting of a single marked point in $T^2$ labelled by the defining representation $V$. We call the first summand the \emph{single skein part} of the skein module, and the second summand the \emph{empty skein part}. The single skein part is considered and its dimension calculated in \S \ref{s-single-skein}. For the empty skein part, we apply the decomposition of \S \ref{s-HH_0-SkAlg}. To obtain explicit dimensions, we need to account for the Weyl group action on the bimodule $A_{\gamma w}$. In \S \ref{s-SL_2-dims-total} we explain how to make these corrections for the $\SL_2$ case and we combine with the single skein dimension to give a formula for the dimension of the whole skein module of $M_\gamma$.

%% file: 4-SL_2-dims/single-skein.tex
In this section we compute the dimension of the space $\HH_0^\gamma(\End_{\SkCat_{\SL_2}(T^2)}(V))$, which is a direct summand of the skein module of $M_\gamma$. We begin by giving a description of the endomorphism algebra. We then give a basis of idempotents for this algebra, and use this to determine the dimension of $\HH_0^\gamma(\End_{\SkCat_{\SL_2}(T^2)}(V))$ in terms of the number of fixed points of $\gamma$ acting on this basis. Here we work with the endomorphism algebra in the Kauffman bracket version of the skein category (see Rmk. \ref{r-SL_2-KBSM-equivalent}).

Let us parameterize $T^2 \times I$ with coordinates $\{ (e^{2\pi i r}, e^{2\pi i s}, t): (r, s, t) \in [0, 1]^3\}$ and assume that all embedded 1-manifolds are given the blackboard framing as in Notation \ref{n-framing}, to allow them to denote skeins.

\begin{prop}
  \label{p-single-skein}
  The map
  \begin{align*}
    \C(q^{1/2})[X, Y]/(X^2 - 1, Y^2 - 1) & \to \End_{\SkCat_{\SL_2}(T^2)}(V)                  \\
    X                                    & \mapsto [\{(e^{2\pi i t}, 0, t) : t \in [0, 1] \}] \\
    Y                                    & \mapsto [\{(0, e^{2\pi i t}, t) : t \in [0, 1]\}]
  \end{align*}
  is an isomorphism of vector spaces.
\end{prop}

\begin{proof}
  We claim there is a surjection $\phi: \C(q^{1/2})[\pi_1(T^2)] \to \End_{\SkCat_{\SL_2}(T^2)}(V)$. This map takes the class of a loop $\alpha : [0, 1] \to T^2$ to the skein given by $\{(x, y, t) : \alpha(t) = (x, y) \}$. Clearly any skein which has just one connected component is in the image of this map, and any other skein can be reduced to a linear combination of such using the skein relations, so the map is a surjection.

  \begin{figure}
    \centering
    \begin{subfigure}[t]{3cm}
      \centering
      \includesvg[width=\textwidth]{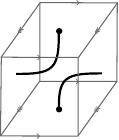}
      \caption{The skein $\phi(X)$, where $\phi$ is the map of Prop. \ref{p-single-skein}.}
      \label{f-X}
    \end{subfigure}
    \hspace{0.5cm}
    \begin{subfigure}[t]{3cm}
      \centering
      \includesvg[width=\textwidth]{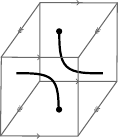}
      \caption{The skein $\phi(X^{-1})$, where $\phi$ is the map of Prop. \ref{p-single-skein}.}
      \label{f-X^{-1}}
    \end{subfigure}
    \hspace{0.5cm}
    \begin{subfigure}[t]{3cm}
      \centering
      \includesvg[width=\textwidth]{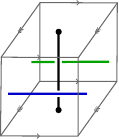}
      \caption{Two skeins $A$ (black and green) and $B$ (black and blue) which are isotopic.}
      \label{f-two-skeins}
    \end{subfigure}
    \caption{Some skeins considered in the proof of Prop. \ref{p-single-skein}.}
  \end{figure}

  Denote by $X$ the class of the usual meridian and by $Y$ the class of the usual longitude. For example, $X$ maps to the skein in Fig. \ref{f-X}. Now, recall that the skein relations can be written in the form
  \begin{equation}
    \label{eq-skein-rel}
    \includesvg[width=10cm]{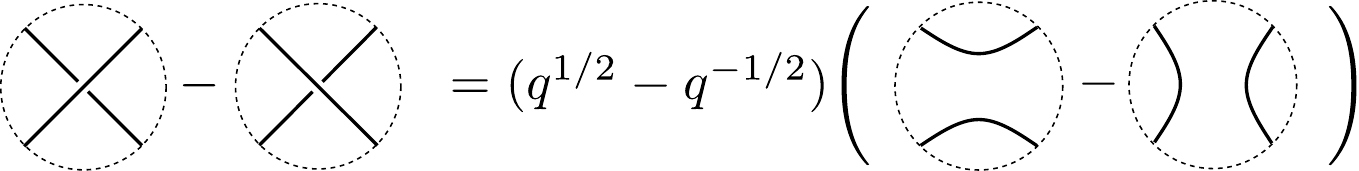}
  \end{equation}
  Then consider the following two skeins:
  \[
    A = \{(e^{2\pi i r}, e^{3\pi i/2}, 1/2) : r \in [0, 1]\} \cup \{(e^{\pi i}, e^{\pi i}, t) : t \in [0, 1]\}
  \]
  and
  \[
    B = \{(e^{2\pi i r}, e^{\pi i/2}, 1/2) : r \in [0, 1]\} \cup \{(e^{\pi i}, e^{\pi i}, t) : t \in [0, 1]\}
  \]
  shown in Fig. \ref{f-two-skeins}. There is an isotopy in $T^2 \times I$ which takes skein $A$ to skein $B$ given by 
  \[
    A_{\theta} = \{(e^{2\pi i r}, e^{3\pi i/2   + \theta}, 1/2) : r \in [0, 1]\} \cup \{(e^{\pi i}, e^{\pi i}, t) : t \in [0, 1]\}, \theta \in [0, \pi]
  \]
  i.e. moving the component $\{(e^{2\pi i r}, e^{3\pi i/2}, 1/2) : r \in [0, 1]\}$ of skein $A$ around the torus in a direction perpendicular to both components of the skein as it is drawn in Fig. \ref{f-two-skeins}. Applying the skein relation (\ref{eq-skein-rel}) we have
  \begin{align*}
    0 = A - B = (q^{1/2} - q^{-1/2})(\phi(X) - \phi(X^{-1}))
  \end{align*}
  so that $\phi(X) = \phi(X^{-1})$, or $\phi(X^2) = 1$. A similar argument shows that $\phi(Y^2) = 1$, and so we have $\ker \phi \supseteq (X^2 - 1, Y^2 - 1)$. In fact we claim this is an equality.

  Since $\ker \phi \supseteq (X^2 - 1, Y^2 - 1)$, we have that $\dim \End_{\SkCat_{\SL_2}(T^2)}(V) \leq 4$. We will argue that $\dim \End_{\SkCat_{\SL_2}(T^2)}(V) \geq 4$, so that $\dim \End_{\SkCat_{\SL_2}(T^2)}(V) = \mathrm{codim} \ker \phi = 4$, from which it follows that $\ker \phi = (X^2 - 1, Y^2 - 1)$ as claimed.

  To give a lower bound on $\dim \End_{\SkCat_{\SL_2}(T^2)}(V)$, we observe that the dimension of an algebra is bounded below by the dimension of its Hochschild homology. So we will give an argument that $\dim \HH_0(\End_{\SkCat_{\SL_2}(T^2)}(V)) = 4$.

  We view $\HH_0(\End_{\SkCat_{\SL_2}(T^2)}(V))$ as a direct summand of $\HH_0(\SkCat_{\SL_2}(T^2)) = \Sk_{\SL_2}(T^3)$. Recall the Kauffman bracket presentation of the skein module: in this case, skeins are equivalence classes of (framed) unoriented links, and so give homology classes in $\H_1(M, \Z/2\Z)$. By the fact that the skein relations relate links with the same $\Z/2\Z$ homology class, as described in \cite{carregaGeneratorsSkeinSpace2017} skein modules are graded by homology with $\Z/2\Z$ coefficients, so that $\HH_0(\SkCat_{\SL_2}(T^2))$ decomposes over $\H_1(T^3, \Z/2\Z) \cong (\Z/2\Z)^3$.

  Moreover, the mapping class group of $T^3$ acts on this space via the surjection $\Z^3 \to (\Z/2\Z)^3$, and it acts transitively on the components which do not correspond to the grading $(0, 0, 0)$, see \cite{carregaGeneratorsSkeinSpace2017}. It has been shown in \cite{gilmerKauffmanBracketSkein2018} that the component in grading $(1, 0, 0)$ is 1-dimensional, from which we see that the subspace graded by $\{(m, n, 1): m, n \in \Z / 2\Z \}$ is 4-dimensional. This is the space of skeins which have nontrivial homology in the direction of the third coordinate, that is, it is $\HH_0(\End_{\SkCat_{\SL_2}(T^2)}(V))$.

  It follows that the dimension of $\End_{\SkCat_{\SL_2}(T^2)}(V)$ is 4, so we have $\ker \phi = (X^2 - 1, Y^2 - 1)$ and $\End_{\SkCat_{\SL_2}(T^2)}(V) \cong \C(q^{1/2})[X, Y]/(X^2 - 1, Y^2 - 1)$.
\end{proof}

Let us write
\begin{align*}
  \C(q^{1/2})[X, Y]/(X^2 - 1, Y^2 - 1) &\cong \C(q^{1/2})[\F_2^2]\\
  X^rY^s  &\mapsto v_{(r,  s)}
\end{align*}
where $v_x$ denotes the $\C(q^{1/2})$-basis element of $\C(q^{1/2})[\F_2^2]$ supported at $x \in \F_2^2$. By Prop. \ref{p-single-skein} we can define a multiplication on this space induced by the multiplication on $\End_{\SkCat_{\SL_2}(T^2)}(V)$: it is given by
\begin{equation}
  \label{eq-SL_2-action-single-skein}
  v_{(r, s)}v_{(t, u)} = v_{(r + t, s + u)}.
\end{equation}
Also from the description of Prop. \ref{p-single-skein} we see that the $\SL_2(\Z)$-action on $\C(q^{1/2})[\F_2^2]$ is given by
\[
  \gamma(v_{(r, s)}) = v_{\bar{\gamma}(r, s)}
\]
where $\bar{\gamma} \in \SL_2(\F_2)$ is the mod 2 reduction of $\gamma$, acting on $\F_2$ in the standard way.

\begin{prop}
  \label{p-single-skein-dims}
  For $\gamma \in \SL_2(\Z)$ and denote by $\bar{\gamma}$ its reduction mod 2. Then we have
  \[
    \dim \HH_0^{\gamma}(\End_{\SkCat_{\SL_2}(T^2)}(V)) = \begin{cases}
      4 & \bar{\gamma} = \bar{\Id}\\
      2 & \tr(\bar{\gamma}) = 0 \mod 2 \text{ and } \bar{\gamma} \neq \bar{\Id}\\
      1 & \tr(\bar{\gamma}) = 1 \mod 2
    \end{cases}.
  \]
\end{prop}

\begin{proof}
  It follows from (\ref{eq-SL_2-action-single-skein}) that imposing the commutator relation $[v_x, v_y]_\gamma$ on $\C(q^{1/2})[\F_2^2]$ is equivalent to imposing
  \[
    (x + y) - (y + \bar{\gamma}(x)) = x - \bar{\gamma}(x)
  \]
  on $\F_2^2$. If $x$ is fixed by $\bar{\gamma}$ then this relation is vacuous, but if $x$ is not fixed then the relation identifies two different basis elements of $\C(q^{1/2})[\F_2^2]$, because $\bar{\gamma}$ is an invertible matrix and so it permutes the elements of $\F_2^2$. Then it follows that the dimension of $\C(q^{1/2})[\F_2^2]$ modulo the $\gamma$-twisted commutator relations is equal to the number of fixed points of $\bar{\gamma}$ acting on $\F_2^2$. The fixed points form a subspace of $\F_2^2$. There are three possibilities:
  \begin{itemize}
    \item There is a 2-dimensional subspace of fixed points if and only if $\bar{\gamma} = \bar{\Id}$. Then there are 4 fixed points.
    \item There is a 1-dimensional subspace of fixed points if and only uf $\bar{\gamma}$ is one of $\begin{pmatrix}1 & 1\\ 0 & 1\end{pmatrix}, \begin{pmatrix}1 & 0\\ 1 & 1\end{pmatrix}$ or $\begin{pmatrix}0 & 1\\ 1 & 0\end{pmatrix}$. That is, when $\tr(\bar{\gamma}) =  0 \mod 2$ and $\bar{\gamma} \neq \bar{\Id}$. Then there are 2 fixed points.
    \item There is a 0-dimensional subspace of fixed points otherwise, i.e. for $\tr(\bar{\gamma}) = 1 \mod 2$. Then there is one fixed point.
  \end{itemize}
\end{proof}

%% file: 4-SL_2-dims/total.tex
We are now ready to give the proof of our main theorem, determining the dimension of the Kauffman bracket skein module for mapping tori of $T^2$.

\begin{thm}
  \label{th-dimension-formula}
  The dimension $\dim \Sk_{\SL_2}(M_\gamma)$ is given by
  \[
    s_\gamma + \frac{\prod_{i = 1}^{r_+}a_i^{+} + 2^{p_{+}}}{2} + \frac{\prod_{i = 1}^{r_-}a_i^{-} + 2^{p_{-}}}{2} 
  \]
  where $p_{\pm} = \# \{a_i^{\pm} \mathrm{ even} : 1 \leq i \leq r_{\pm} \}$ for $a_i^\pm$ the invariant factors of $\Id \mp \gamma : \Lambda^2 \to \Lambda^2$ and $r_{\pm}$ the rank of this map, and 
  \[
    s_\gamma = \begin{cases}
      4 & \bar{\gamma} = \bar{\Id}\\
      2 & \tr(\bar{\gamma}) = 0 \mod 2 \text{ and } \bar{\gamma} \neq \bar{\Id}\\
      1 & \tr(\bar{\gamma}) = 1 \mod 2
    \end{cases}
  \]
 for $\bar{\gamma}$ the matrix for $\gamma$ reduced modulo 2 in each entry.
\end{thm}

\begin{proof}
  Firstly, by Cor. \ref{c-skmod-decomp} with $G = \SL_2, W = \Z /2\Z$ we see that the dimension is given by the sum $s_\gamma + d_{+} + d_{-}$
  where $s_\gamma$ is the dimension of the single skein part, which was computed in Prop. \ref{p-single-skein-dims}, and $d_{\pm}  = \dim \HH_0(A, A_{\pm \gamma})_{\Stab(\pm 1)}$. It therefore suffices to understand $d_{\pm}$. Before taking coinvariants, we have that
  \begin{equation*}
    \HH_0(A, A_{\pm \gamma}) \cong \C(q^{1/2})\left[ \coker(\Id \mp \gamma)_{\mathrm{tors}} \right] \cong \C(q^{1/2})\left[ \bigoplus_{i = 1}^{r_{\pm}} \Z / a_i^{\pm} \Z \right]
  \end{equation*}
  where the first isomorphism is Cor. \ref{c-dim-HH_0-before-stab} and second isomorphism is standard. From this it follows that the dimension before taking coinvariants of each component $\HH_0(A, A_{\pm\gamma})$ is given by $\prod_{i = 1}^{r_\pm} a_i^{\pm}$. All that remains is to account for the action of the centralizers.

  We note that $\Stab(w) = W$ for all $w \in W$ since $W = \Z/2\Z$, and observe that $W$ acts on each component of the lattice $\Lambda^2$ by negation. We need to count the number of orbits of the induced action on the set $\bigoplus_{i = 1}^{r_{\pm}} \Z / a_i^{\pm} \Z$, which we consider as a rectangular subset of $\Lambda^2$. Since the orbits will be generically of size 2 then this will be approximately $\frac{1}{2}(\prod_{i = 1}^{r_\pm}a_i^{\pm})$. The precise number of orbits will depend on the number of fixed points of the action of $W$.

  Observe that 0 is always fixed. If none of the $a_i^{\pm}$ are even, then this is the only fixed point. If precisely one of the $a_i^{\pm}$ is even then this forces there to be another fixed point on the $i$-axis of the set. If two of the $a_i^{\pm}$ are even then there are four fixed points, at the four corners of the set. Then the precise number of orbits is given by
  \[
    \frac{\prod_{i = 1}^{r_\pm}a_i^{\pm} + 2^{p_\pm}}{2}
  \]
  where $2^{p_\pm}$ counts the number of fixed points, for $p_{\pm}$ as given in the statement of the theorem. 
\end{proof}

%% file: A-Appendix/code.tex
Here we include Sage code for implementing the formula of Thm. \ref{th-dimension-formula}. This implementation can be found together with some precomputed dimensions at the following repository:\\\url{https://github.com/PatrickKinnear/skein-dimensions.git}.

\begin{lstlisting}[float=!h]
def get_dim_single_skein(gamma):
    I = matrix(ZZ, 2, [1, 0, 0, 1])
    P = matrix(ZZ, 2, [1, 1, 1, 0])
    Q = matrix(ZZ, 2, [0, 1, 1, 1])

    if gamma % 2 == I:
        return 4
    elif gamma % 2 == P or gamma % 2 == Q:
        return 1
    else:
        return 2

def get_dim_empty_skein(gamma):
    I = matrix(ZZ, 2, [1, 0, 0, 1])

    D_plus, U_plus, V_plus = (I - gamma).smith_form()
    a_plus = [a for a in D_plus.diagonal() if a != 0]

    D_minus, U_minus, V_minus = (I + gamma).smith_form()
    a_minus = [a for a in D_minus.diagonal() if a != 0]

    p_plus = len([a for a in a_plus if a%2 == 0])
    p_minus = len([a for a in a_minus if a%2 == 0])

    return [(prod(a_plus) + 2**(p_plus))/2, \
            (prod(a_minus) + 2**p_minus)/2]

def skein_dimension(gamma):
    s_gamma = get_dim_single_skein(gamma)
    d_plus, d_minus = get_dim_empty_skein(gamma)

    return s_gamma + d_plus + d_minus
\end{lstlisting}

%% file: A-Appendix/table.tex
\begin{table}[H]
  \begin{minipage}[t]{0.49\textwidth}
    \begin{center}
      \begin{tabular}[t]{l || r r r | r}
            $\gamma$ & $s_\gamma$ & $d_{+}$ & $d_{-}$ & Total\\
            \hline
            $R^{	1	}L^{	1	}$ &	1	&	1	&	3	&	5	\\
            $R^{	1	}L^{	2	}$ &	2	&	2	&	4	&	8	\\
            $R^{	1	}L^{	3	}$ &	1	&	2	&	4	&	7	\\
            $R^{	1	}L^{	4	}$ &	2	&	3	&	5	&	10	\\
            $R^{	1	}L^{	5	}$ &	1	&	3	&	5	&	9	\\
            $R^{	1	}L^{	6	}$ &	2	&	4	&	6	&	12	\\
            $R^{	1	}L^{	7	}$ &	1	&	4	&	6	&	11	\\
            $R^{	1	}L^{	8	}$ &	2	&	5	&	7	&	14	\\
            $R^{	1	}L^{	9	}$ &	1	&	5	&	7	&	13	\\
            $R^{	1	}L^{	10	}$ &	2	&	6	&	8	&	16	\\
            $R^{	2	}L^{	2	}$ &	4	&	4	&	6	&	14	\\
            $R^{	2	}L^{	3	}$ &	2	&	4	&	6	&	12	\\
            $R^{	2	}L^{	4	}$ &	4	&	6	&	8	&	18	\\
            $R^{	2	}L^{	5	}$ &	2	&	6	&	8	&	16	\\
            $R^{	2	}L^{	6	}$ &	4	&	8	&	10	&	22	\\
            $R^{	2	}L^{	7	}$ &	2	&	8	&	10	&	20	\\
            $R^{	2	}L^{	8	}$ &	4	&	10	&	12	&	26	\\
            $R^{	2	}L^{	9	}$ &	2	&	10	&	12	&	24	\\
            $R^{	2	}L^{	10	}$ &	4	&	12	&	14	&	30	\\
            $R^{	3	}L^{	3	}$ &	1	&	5	&	7	&	13	\\
            $R^{	3	}L^{	4	}$ &	2	&	7	&	9	&	18	\\
            $R^{	3	}L^{	5	}$ &	1	&	8	&	10	&	19	\\
            $R^{	3	}L^{	6	}$ &	2	&	10	&	12	&	24	\\
            $R^{	3	}L^{	7	}$ &	1	&	11	&	13	&	25	\\
            $R^{	3	}L^{	8	}$ &	2	&	13	&	15	&	30	\\
            $R^{	3	}L^{	9	}$ &	1	&	14	&	16	&	31	\\
            $R^{	3	}L^{	10	}$ &	2	&	16	&	18	&	36	\\
            $R^{	4	}L^{	4	}$ &	4	&	10	&	12	&	26
        \end{tabular}
    \end{center}
    \end{minipage}
    \hfill
    \begin{minipage}[t]{0.49\textwidth}
      \begin{center}
        \begin{tabular}[t]{l || r r r | r}
            $\gamma$ & $s_\gamma$ & $d_{+}$ & $d_{-}$ & Total\\
            \hline
            $R^{	4	}L^{	5	}$ &	2	&	11	&	13	&	26	\\
            $R^{	4	}L^{	6	}$ &	4	&	14	&	16	&	34	\\
            $R^{	4	}L^{	7	}$ &	2	&	15	&	17	&	34	\\
            $R^{	4	}L^{	8	}$ &	4	&	18	&	20	&	42	\\
            $R^{	4	}L^{	9	}$ &	2	&	19	&	21	&	42	\\
            $R^{	4	}L^{	10	}$ &	4	&	22	&	24	&	50	\\
            $R^{	5	}L^{	5	}$ &	1	&	13	&	15	&	29	\\
            $R^{	5	}L^{	6	}$ &	2	&	16	&	18	&	36	\\
            $R^{	5	}L^{	7	}$ &	1	&	18	&	20	&	39	\\
            $R^{	5	}L^{	8	}$ &	2	&	21	&	23	&	46	\\
            $R^{	5	}L^{	9	}$ &	1	&	23	&	25	&	49	\\
            $R^{	5	}L^{	10	}$ &	2	&	26	&	28	&	56	\\
            $R^{	6	}L^{	6	}$ &	4	&	20	&	22	&	46	\\
            $R^{	6	}L^{	7	}$ &	2	&	22	&	24	&	48	\\
            $R^{	6	}L^{	8	}$ &	4	&	26	&	28	&	58	\\
            $R^{	6	}L^{	9	}$ &	2	&	28	&	30	&	60	\\
            $R^{	6	}L^{	10	}$ &	4	&	32	&	34	&	70	\\
            $R^{	7	}L^{	7	}$ &	1	&	25	&	27	&	53	\\
            $R^{	7	}L^{	8	}$ &	2	&	29	&	31	&	62	\\
            $R^{	7	}L^{	9	}$ &	1	&	32	&	34	&	67	\\
            $R^{	7	}L^{	10	}$ &	2	&	36	&	38	&	76	\\
            $R^{	8	}L^{	8	}$ &	4	&	34	&	36	&	74	\\
            $R^{	8	}L^{	9	}$ &	2	&	37	&	39	&	78	\\
            $R^{	8	}L^{	10	}$ &	4	&	42	&	44	&	90	\\
            $R^{	9	}L^{	9	}$ &	1	&	41	&	43	&	85	\\
            $R^{	9	}L^{	10	}$ &	2	&	46	&	48	&	96	\\
            $R^{	10	}L^{	10	}$ &	4	&	52	&	54	&	110  \\
        \end{tabular}
    \end{center}
  \end{minipage}
  \begin{center}
  \caption{Kauffman bracket skein module dimensions for hyperbolic mapping classes corresponding to sequences of length 2.}
  \label{t-length-2-dims}
  \end{center}
\end{table}

%% file: A-Appendix/derivation.tex
We recall the following elementary facts about the invariant factors of an integer matrix $P$. The invariant factors satisfy $a_i | a_{i+1}$ and can be calculated by $a_i = d_i/d_{i-1}$ where $d_i$ is the \emph{$i$-th determinant divisor}, the greatest common divisor of the determinants of all $i \times i$ minors of $P$ (we take the convention that $d_0 = 1$). Moreover if $P$ is full rank, then $\prod_i a_i = d_r = |\det(P)|$. Using these facts, we can show the following.

\begin{cor}
    For $\gamma = \begin{pmatrix} a & b\\ c & d\end{pmatrix} \in \SL_2(\Z)$ with $|\tr(\gamma)| > 2$, we have 
    \[
    \dim \Sk_{\SL_2}(M_{\gamma}) = |\tr(\gamma)| + 2^{c(\gamma) + 1}
  \]
  where
  \[
    c(\gamma) = \#\{m \in \{\mathrm{gcd}(a-1, b, c, d-1), \tr(\gamma)\} : m \text{ even} \}.
  \]
\end{cor}

\begin{proof}
We see that $\Id \mp \gamma$ has full rank since
\[
  \det(\Id \mp \gamma) = \chi_{\gamma}(\pm 1) = 2 \mp \tr(\gamma) \neq 0.
\]

We write $a_i^{\pm}, d_i^{\pm}$ for the invariant factors and determinant divisors respectively of $\Id \mp \gamma$. We have
\begin{align*}
  \frac{1}{2}(\prod a_i^{+} + \prod a_i^{-}) &= \frac{1}{2}(|\det(\Id - \gamma)| + |\det(\Id + \gamma)|)\\ 
  &= \frac{1}{2}(|2 -  \tr(\gamma)| + |2 + \tr(\gamma)|)\\
  &= |\tr(\gamma)|.
\end{align*}

Let us consider the numbers $p_{\pm}$ from Theorem \ref{th-dimension-formula}. Writing $a_1^{\pm}$ as the greatest common divisor of the matrix entries, we notice that $2 | a_1^{+} \iff 2 | a_1^{-}$. Then since  $a_1^{\pm} | a_2^{\pm}$, we see $p_{+} = 2 \iff p_{-} = 2$.

Assuming $2 \nmid a_1^{\pm}$, it follows that $2 | a_2^{\pm} \iff 2 | d_2^{\pm}$. In this case, writing $d_2^{\pm} = |2 \mp \tr(\gamma)|$, we see $2 | a_2^{+} \iff 2 | \tr(\gamma) \iff 2 | a_2^{-}$, and so $p_{+} = 1 \iff p_{-} = 1$. It then follows that $p_{+} = 0 \iff p_{-} = 0$. 

Then we see that 
\[
  2^{p_{+} - 1} + 2^{p_{-} - 1} = 2^{c(\gamma)}
\]
for $c(\gamma)$ as above. It is also easy to see that $s_{\gamma} = 2^{c(\gamma)}$. Then it follows that
\begin{align*}
  s_{\gamma} + \frac{\prod_{i = 1}^{r_+}a_i^{+} + 2^{p_{+}}}{2} + \frac{\prod_{i = 1}^{r_-}a_i^{-} + 2^{p_{-}}}{2} &= 2^{c(\gamma)} + 2^{p_+ - 1} + 2^{p_1 - 1} + \frac{1}{2}(\prod_{i = 1}^{2} a_i^{+} + \prod_{i = 1}^{2} a_i^{-})\\
  &= 2^{c(\gamma)} + 2^{c(\gamma)} + |\tr(\gamma)|\\
  &= 2^{c(\gamma) + 1} + |\tr(\gamma)|
\end{align*}
which verifies the formula.
\end{proof}